\documentclass{jfm}

\usepackage{graphicx}
\usepackage{natbib}

\ifCUPmtlplainloaded \else
  \checkfont{eurm10}
  \iffontfound
    \IfFileExists{upmath.sty}
      {\typeout{^^JFound AMS Euler Roman fonts on the system,
                   using the 'upmath' package.^^J}
       \usepackage{upmath}}
      {\typeout{^^JFound AMS Euler Roman fonts on the system, but you
                   dont seem to have the}%
       \typeout{'upmath' package installed. JFM.cls can take advantage
                 of these fonts,^^Jif you use 'upmath' package.^^J}%
       \providecommand\upi{\pi}%
      }
  \else
    \providecommand\upi{\pi}%
  \fi
\fi

\ifCUPmtlplainloaded \else
  \checkfont{msam10}
  \iffontfound
    \IfFileExists{amssymb.sty}
      {\typeout{^^JFound AMS Symbol fonts on the system, using the
                'amssymb' package.^^J}%
       \usepackage{amssymb}%
         \let\leq=\leqslant
         \let\geq=\geqslant
      }{}
  \fi
\fi

\ifCUPmtlplainloaded \else
  \IfFileExists{amsbsy.sty}
    {\typeout{^^JFound the 'amsbsy' package on the system, using it.^^J}%
     \usepackage{amsbsy}}
    {\providecommand\boldsymbol[1]{\mbox{\boldmath $##1$}}}
\fi

\usepackage{float}           
\usepackage{todonotes}       
\usepackage{cancel}          
\usepackage{cases}           

\usepackage{subfig}          
\usepackage{textcomp}	     
\usepackage{multirow}	     

\providecommand\bnabla{\boldsymbol{\nabla}}
\providecommand\bcdot{\boldsymbol{\cdot}}

\newsavebox{\astrutbox}
\sbox{\astrutbox}{\rule[-5pt]{0pt}{20pt}}

\newcommand\p{\ensuremath{\partial}}

\newcommand\etal{\mbox{\textit{et al.}}}

\newcommand\eg{e.g.\ }
\newcommand\ie{i.e.\ }
\newcommand\ci{\mathrm{i}}

\newcommand{\ds}{\displaystyle}
\newcommand{\norm}[1]{|#1|}

\newcommand{\bc}{{\boldsymbol{c}}}
\newcommand{\be}{{\boldsymbol{e}}}

\newcommand{\bn}{{\boldsymbol{n}}}

\newcommand{\bq}{{\boldsymbol{q}}}

\newcommand{\bx}{{{\boldsymbol{x}}}}
\newcommand{\by}{{{\boldsymbol{y}}}}

\newcommand{\bv}{{\boldsymbol{v}}}

\newcommand{\bX}{ {\boldsymbol{X}}}
\newcommand{\bY}{ {\boldsymbol{Y}}}

\newcommand{\bbeta}{ \boldsymbol{\beta}}
\newcommand{\bxi}{ \boldsymbol{\xi}}
\newcommand{\btau}{ \boldsymbol{\tau}}

\newcommand{\bPi}{\boldsymbol{\Pi}}
\newcommand{\bmu}{\boldsymbol{\mu}}

\newcommand{\sA}{\mathsf{A}}

\newcommand{\sD}{\mathsf{D}}
\newcommand{\sH}{\mathsf{H}}
\newcommand{\sJ}{\mathsf{J}}

\newcommand{\sP}{\mathsf{P}}

\newcommand{\sU}{\mathsf{U}}

\newcommand{\sW}{\mathsf{W}}

\newcommand{\absxxi}{\vert \boldsymbol{x}-\boldsymbol{\xi} \vert}

\newcommand{\partialtroispm}{\partial_{x_3}p_\pm}

\newcommand{\partialtrois}{\partial_{x_3}p}

\newcommand{\hsD}{\widehat{\sD}}

\newcommand{\hsDpm}{\widehat{\sD}_\pm}

\newcommand{\homega}{\widehat{\Omega}}
\newcommand{\hsigma}{\widehat{\Sigma}}
\newcommand{\hKR}{\widehat{K_R}}
\newcommand{\hx}{\hat{x}}
\newcommand{\bhx}{{{\boldsymbol{\hx}}}}

\newcommand{\pizero}{\Pi^{(0)}}

\newcommand{\piun}{\Pi^{(1)}}

\newcommand{\dif}{\mathrm{d}}
\newcommand{\epss}{\varepsilon}

\newenvironment{systeme}{\begin{equation}
\left\{ \begin{array}{ll} \displaystyle }{\end{array} \right.
\end{equation}}

\renewcommand{\Re}{\mbox{Re}}
\renewcommand{\Im}{\mbox{Im}}
\renewcommand{\circ}{\textit{o}}
\renewcommand{\bigcirc}{\textit{O}}

\newtheorem{lemma}{Lemma}

\title[]{Effective conditions for the reflection of an acoustic wave by low-porosity perforated plates}

\author[S. Laurens, E. Piot, A. Bendali, M'B. Fares and S. Tordeux ]%
{S.\ns L\ls A\ls U\ls R\ls E\ls N\ls S$^{1,3}$%
  \thanks{Email address for correspondence: sophie.laurens@cerfacs.fr},\ns
E.\ns P\ls I\ls O\ls T$^2$, A. \ns B\ls E\ls N\ls D\ls A\ls L\ls I$^{1,3}$,
M'B. \ns F\ls A\ls R\ls E\ls S$^3$,\ns
\and S. \ns T\ls O\ls R\ls D\ls E\ls U\ls X$^4$}

\affiliation{$^1$IMT, University of Toulouse, INSA, 135 avenue de Rangueil, F-31077, Toulouse, France\\[\affilskip]
$^2$ONERA - The French Aerospace Lab, F-31055, Toulouse, France\\[\affilskip]
$^3$CERFACS, 42 avenue Gaspard Coriolis, F-31100, Toulouse, France\\[\affilskip]
$^4$ INRIA \& University of Pau, LMA, avenue de l'Universit\'e, F-64000, Pau, France}

\pubyear{2012}
\volume{650}
\pagerange{119--126}
\date{?; revised ?; accepted ?. - To be entered by editorial office}
\begin{document}

\maketitle

\begin{abstract}
This paper describes an investigation of the acoustic properties of a
rigid plate with a periodic pattern of holes, in a compressible,
ideal, inviscid fluid in the absence of mean flow. Leppington and
Levine (\textit{Reflexion and transmission at a plane screen with
  periodically arranged circular or elliptical apertures}, J. Fluid
Mech., 1973, p.109-127) obtained an approximation of the reflection and transmission coefficients
of a plane wave incident on an infinitely thin plate with a
rectangular array of perforations, assuming that a characteristic
size of the perforations is negligible relative to that of the unit
cell of the grating, itself assumed to be negligible relative to the
wavelength. One part of the present
study is of methodological interest. It establishes that it is possible to
extend their approach to thick plates with a
skew grating of perforations, thus confirming recent results in Bendali \textit{et al.}
(\textit{Mathematical justification of the Rayleigh conductivity model
  for perforated plates in acoustics}, SIAM J. Appl. Math., 2013), but
in a much simpler way without using complex matched asymptotic expansions of
the full wave or to a grating of multipoles. As is well-known, effective compliances for the plate can then be derived from the corresponding approximations of the reflection and transmission
 coefficients. These compliances are expressed in terms of the
Rayleigh conductivity of an isolated perforation. Consequently, in one other part
of the present study, the methodology recently introduced in  Laurens
\textit{et al.} (\textit{Lower and upper bounds for the Rayleigh
  conductivity of a perforated plate}, ESAIM:M2AN, 2013) to obtain
sharp bounds for the Rayleigh conductivity has been extended to include the case
for which the openings of the perforations on the upper and lower sides of
the plate are elliptical in shape. This not only enables the determination of
these bounds and of the associated reflection and transmission coefficients for actual plates with tilted perforations but also
yields single expressions covering all usual cases of perforations:
straight or tilted with a circular or an elliptical cross-section.
\end{abstract}

\begin{keywords}
\end{keywords}

\section{Introduction}

Perforated plates and screens are widely used in engineering systems due to their ability to absorb sound or to reduce sound transmission, in a variety of applications
 including room acoustics and aeroacoustics \citep{1994_Ingard,1992_book_beranek}. They can be used as protective layers of porous materials
   \citep{1993_book_allard} to form sandwich structures in aircraft fuselages \citep{1999_JASA_Abrahams} or as facing layers of liners.
    In this case, the perforated plates are either backed by honeycomb cells which are mounted on a rigid backplate \citep{1998_JASA_Maa} or by a plenum 
    acting as a resonant cavity in gas turbine combustion chambers \citep{1999_book_Lefebvre}. 

The reflection and transmission properties of such panels perforated with a periodic array of holes can be modeled directly, as was done by \citet{1973_JFM_Leppington} 
by  means of an integral equation
 or later by the method of matched asymptotic expansions (see
 \citet{bendali2012} and \citet{1990_PRSL_Leppington} for the case of
 elastic plates). The acoustic effect of perforated plates
 is usually addressed by studying the acoustic behaviour of a single aperture in an infinite wall and by designing a homogeneous model for the whole plate,
  mainly by means of an averaging procedure. In this way, the Leppington and
  Levine approach was later extended by \citet{1980_JFM_Howe} and \citet{1992_JSV_Dowling} to include a mean tangential or bias flow to the perforated screen.
In these papers, the asymptotic solution was sought for the long wavelength limit and under the low-porosity assumption, i.e., 
the size of each perforation is assumed to be small in comparison to the
 distance between two neighbouring apertures of the periodic array,
 viewed in turn as being much smaller than the wavelength of the
 impinging sound waves. In \cite{bendali2012}, it is the whole wave which is expanded by means of the method
of matched asymptotic expansions. As a result, the
procedure for determining the two top-order terms of the asymptotic
expansions of the reflection and transmission coefficients relative
to the small parameter, characterizing the ratio of a characteristic
size of the perforation to the spacing between two successive
perforations, is particularly complex.

All these homogenisation processes rely on a model of the acoustic behaviour of the apertures, which depends on the actual problem being considered.
In the studies for aeronautical turbofan engines \citep[see for example][for a recent review]{tam_jsv_08}, acoustic impedance is the relevant parameter, 
more precisely,  the ratio of the fluctuating pressure  across the hole to the normal fluctuating velocity through it. The real part of the impedance is 
denoted as a resistance and its imaginary part as a reactance. On the contrary, with respect to gas turbines
 \citep[see][for instance]{2011_ASME_Andreini,2012_JSV_Scarpato}, the acoustic behaviour of the aperture is described in terms of  its compliance, termed in this context
  `Rayleigh conductivity', which characterizes the fluctuating volume flow rate through the hole as a function of the difference in unsteady pressure between each side
   of the plate.

The physical phenomenon involved in the acoustic behaviour of an aperture can be broken down into three mechanisms. The first one is purely inviscid and is associated
 with the sound radiation of the perforation and the distortion of the acoustic flow at the plate surface. This flow can be assumed to be potential and governed by 
 the Helmholtz equation, under the assumption that a characteristic size  of the hole is small in comparison to the acoustic wavelength. Consequently,
  no acoustic damping can be predicted by this approach, which means that the impedance $Z$ of the aperture is purely imaginary, while its Rayleigh conductivity $K_R$
   is purely real. This mechanism is generally explained by considering that the flow within the hole behaves as a small piston of air, whose length is 
   greater than
    the depth of the hole because of inertial effects. This amounts to increasing the mass of the vibrating air and is accounted for by using correction lengths, which
     need to be added to the plate thickness. For a circular aperture of radius $r$, \citet{rayleigh} gave a lower and an upper bound for this end correction:
       $\upi r/2\approx 1.57r$ and $16r/3\upi\approx 1.70r$. Later, \citet{1969_JSV_Morfey} derived approximated expressions for openings of arbitrary shape at low
        sound frequencies.  
The second mechanism is linked to viscous effects occurring within the boundary layers which develop at the perforation's inner walls,
 and around the perforation edges at the plate surface. This induces a real part for the acoustic impedance (i.e., an imaginary part of the Rayleigh conductivity),
  which is linked to an absorption of acoustic energy. Viscous effects also modify the imaginary part of the impedance, which is usually modelled by viscous end corrections
   \citep{Melling_JSV,1998_JASA_Maa}.
Finally, the third mechanism consists of vortex shedding, which converts acoustical into mechanical energy,  subsequently dissipated into heat. This may be due to
 nonlinear effects, when high-amplitude sound waves impinge on the aperture in the absence of mean flow \citep{Cummings_84}, or to the effect of mean flow.
  In the presence of a grazing mean flow, empirical \citep{1975_JSV_Guess, 1998_JSV_Kirby}, semi-empirical \citep{1986_Acustica_Cummings, 1987_JSV_Cummings} or
   analytical \citep{1996_PRSL_Howe} models are used to find the impedance or the Rayleigh conductivity of the aperture. It is shown that the resistance induced
    by a high-speed grazing flow is much larger than the viscous resistance obtained without flow, consequently an increased loss of acoustic power is generally observed.
However the phenomena involved are quite complex and have not been accurately and comprehensively modelled. 
On the contrary, \citet{Howe_79} developed an analytical model for the aperture Rayleigh conductivity which proved to be quite
  reliable \citep{1990_JFM_Hughes,jing_jasa_99,2003_JFM_Eldredge,2004_AIAAJ_Bellucci,2009_JCP_Mendez,2012_JSV_Scarpato} for cases in which the fluid flows
  through the perforation (i.e. a bias flow),
 even when mechanisms similar to the grazing case are involved. It consists in multiplying the no-flow Rayleigh
   conductivity of a circular aperture in a zero-thickness plate (taken as $K_R=2r$, i.e., an inertial end correction of $\upi r/2$) by a complex function of the mean bias flow 
   Mach number. The Reynolds number within the perforation is assumed to be large enough so that viscous effects can be neglected except near the edges, where flow separation
    takes place. Later, this model was improved to take into account the plate thickness, by adding it to the inertial end correction \citep{2000_AIAAJ_Jing} or
     by modifying the flow-dependent function \citep{2005_JFS_Howe}. However these expressions only deal with the simplified geometry of a non-tilted circular aperture, while in actual combustors the apertures are tilted downstream and their cross-section in the plate plane is elliptical in shape.
      Recent papers have numerically \citep{2007_AIAA_Eldredge} or experimentally \citep{2011_ASME_Andreini} investigated such a realistic geometry, by using Howe's model with a modified thickness and by neglecting some other geometric effects. 
      Finally, only a few studies have addressed the acoustic response of a perforation in the complex configuration in which an exiting bias flow interacts with a grazing flow, with the exception of \citet{2002_JSV_Sun},  who attempted  to correlate experimental results with an empirical model.

If the low-porosity assumption (aperture dimension $\ll$ spacing of the array) is verified and the long wavelength limit is considered (spacing of the array $\ll$ wavelength), the  device behaves acoustically as though each perforation were isolated. 
In this case, the Rayleigh conductivity (or impedance) of a single aperture, as defined above, is used to define an expression for the compliance (or impedance) of a
 homogeneous screen of such apertures, by assuming that the fluctuating volume flux through each hole is uniformly distributed over each cell of the array
  of perforations \citep[see][]{2003_JFM_Eldredge, 1990_JFM_Hughes,Melling_JSV, 1998_JASA_Maa, 2007_JSV_Atalla}. The layout of the array of holes is not involved in any way.
When the porosity $\sigma$ of the perforated plate is greater than about $4\%$, interaction between the apertures can no longer be neglected. The problem of acoustic 
interaction was solved by \citet{1941_Fok} for the case of an infinitely thin plate with circular perforations. He introduced a porosity-depending interaction function $\psi$ and 
he derived a modified Rayleigh conductivity defined as $K_R'=2r\psi(\sqrt{\sigma})$. \citet{Melling_JSV} used Fok's function $\psi$
 for modifying the inertial end correction (taken equal to $16r/3\upi$ in his model) while other authors \citep{1975_JSV_Guess,1992_book_beranek,1993_book_allard} preferred to use empirical expressions. Recently, interaction effects were studied numerically and experimentally by \citet{2007_JSV_Lee} for the case of perforated plates with a bias flow.

In this paper, we focus on a few questions arising from the previous discussion, which showed that all the impedance or conductivity models are based on the
 inertial end corrections given by \citet{rayleigh} for a single circular aperture in an infinitely thin plate. It should be noted that in room acoustics and liner
  applications, it is the $16r/3\upi$ bound which is used, probably introduced by \citet{1953_JASA_Ingard}. On the contrary, Howe's model \citep{Howe_79} with an end 
  correction of $\upi r/2$ seems to be the one most widely used in papers dealing with gas turbines. Recently, \citet{laurens2012} developed a rigorous mathematical
   framework in order to derive accurate bounds of the Rayleigh conductivity for various hole geometries, but these results do not cover the case 
   of perforated 
   plates whose perforations have circular-cross sections and are tilted, that are in common use in gas turbines. 
 The additional difficulty comes from the openings of the perforations which are then elliptical in shape. Therefore, in the first sections of the present paper, we derive  
 an analytical formula expressing the bounds of the Rayleigh conductivity for almost all the situations which can be encountered in practice: straight or tilted 
 perforation, with either circular or elliptical cross-sections. The associated end corrections are also given. Then, in section \ref{section:reseau}, we show that the two 
 expansions of first and second order of the reflection and transmission coefficients, derived in a fairly complex way by \citet{bendali2012}, can be obtained
  in a much simpler manner by
an extension of the \cite{1973_JFM_Leppington} approach. This method is, in the opinion of the authors, more suitable for a fluid mechanics context 
and it can be easily adapted to other models of flows at the level of the perforations. Moreover, as these coefficients are expressed in terms of the Rayleigh conductivity of a single perforation, 
the effective reflection coefficient of perforated plates commonly used in aeronautics can thus be derived.

\section{Description of the perforated plate and governing equations}\label{section:notations}

\begin{figure}
    \centering   
 \subfloat[Incident and scattered acoustic waves]{\includegraphics[width = 0.475\linewidth]{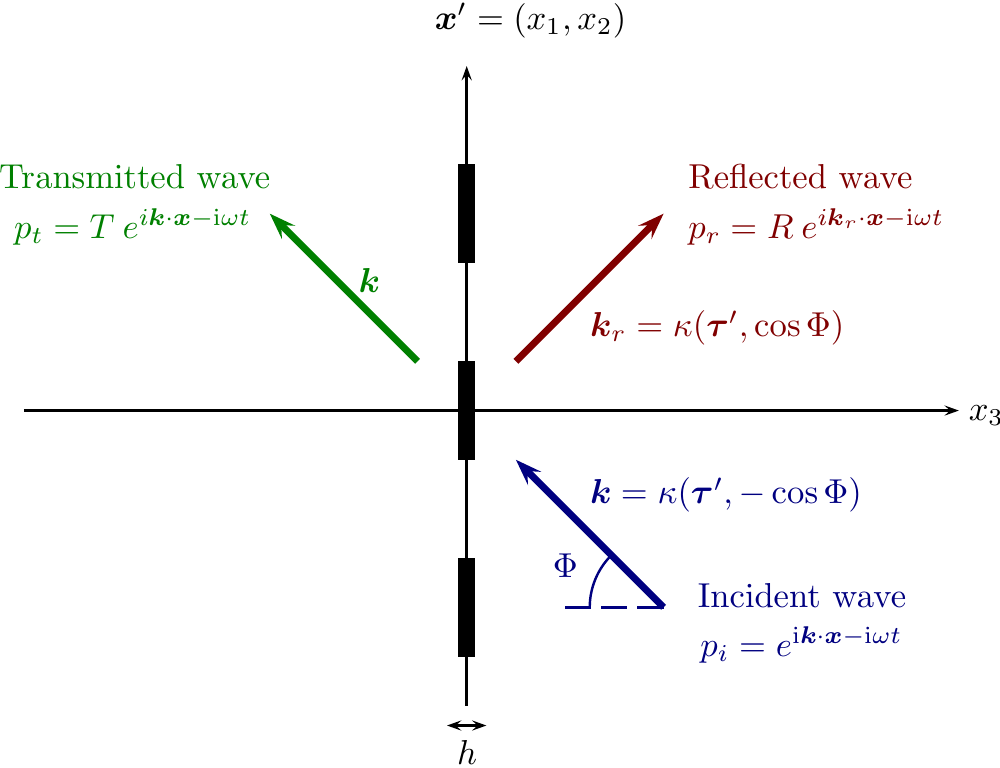}  \label{fig:plaque-reseau}} \hspace{0.75em}
\subfloat[Sketch showing layout of the holes in the lattice]{\includegraphics[width = 0.475\linewidth]{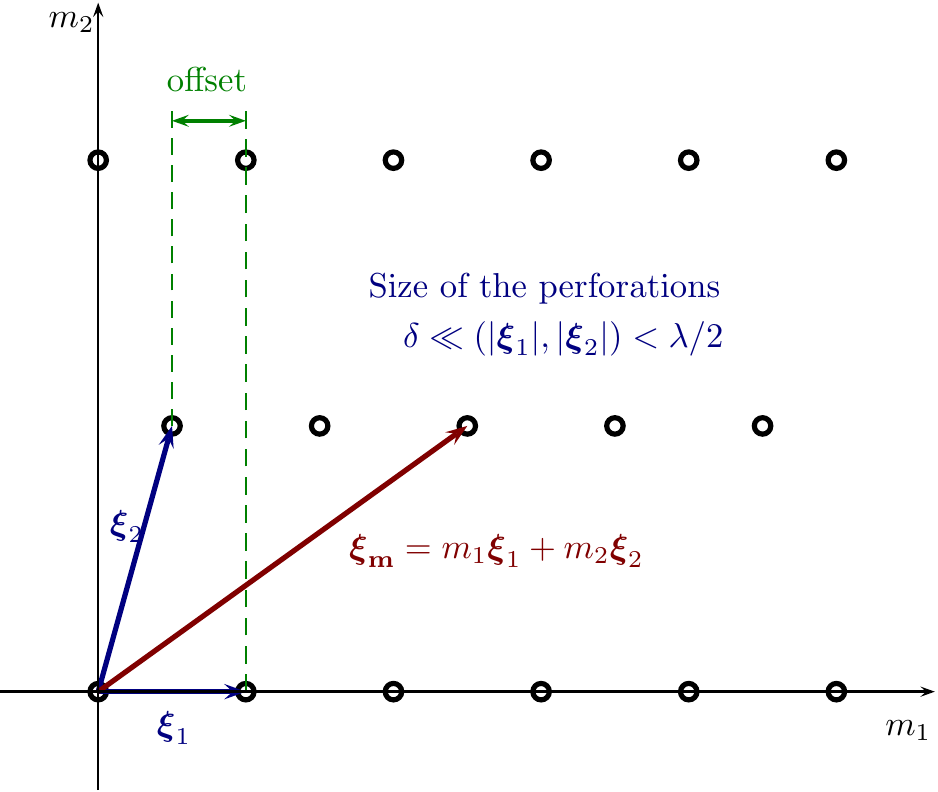} \label{fig:plaque-schema}
}
    \caption{The geometry of the perforated screen}
 \label{fig:plaque}
\end{figure}

Consider the situation shown in figure \ref{fig:plaque} of a perforated screen in a compressible ideal inviscid fluid in the absence of mean flow.
A time-periodic acoustic plane wave is obliquely incident to an infinitely rigid plate of thickness $h$ that lies in the plane $x_3=0$ (see figure \ref{fig:plaque-reseau}) at an angle $\Phi$ to the normal direction to the plate, with $-\pi/2 < \Phi < \pi/2$. Throughout the paper, the coordinate system $\bx$ is decomposed into the two-dimensional variable $\bx' =(x_1, x_2)$ and the coordinate $x_3$ along the normal to the plate. The pressure variation of the incident wave is given by
\begin{equation}\label{eq:pi}
p_i =\exp \left( \ci \kappa (\btau' \bcdot \bx' - x_3 \cos{\Phi} )-\ci \omega t \right)
\end{equation}
with $\btau=(\btau',-\cos{\Phi})$ the cosine directors of the direction of propagation. The wavenumber $\kappa=\omega/c_0$ expresses the ratio of the angular frequency $\omega$ to the speed of sound $c_0$. At a great enough distance from the perforated plate, the scattered pressure field can be decomposed into a reflected wave 
\begin{equation}\label{eq:pr}
p_r =R \:\exp \left( \ci \kappa (\btau' \bcdot \bx' + x_3 \cos{\Phi} )-\ci \omega t \right)
\end{equation}
and a transmitted one
\begin{equation}\label{eq:pt}
p_t =T\: \exp \left( \ci \kappa (\btau' \bcdot \bx' - x_3 \cos{\Phi} )-\ci \omega t \right)
\end{equation}
where $R$ and $T$ are complex constants, which depend on the acoustic properties of the perforated  plate.
The time multiplicative factor $e^{-\ci \omega t}$ is hereafter suppressed by linearity.

In the fluid domain, the pressure field satisfies the homogeneous Helmholtz equation and the cancellation of normal velocity on the rigid walls, equivalently expressed from the linearized momentum equation, by cancelling the pressure normal derivative. Thus, the acoustic wave is governed by the following boundary-value problem:
 \begin{equation} \label{pb:Helmholtz}
\left\{
\begin{array}{ll}
\displaystyle
  \nabla^2 p +\kappa^2 p = 0, & \mbox{in the fluid domain }, \\
   \p_\bn p = 0, & \mbox{on the rigid walls}, \\
 p =  p_i + p_r, & \mbox{as\ }  x_3 \rightarrow +\infty, \\
 p =  p_t, & \mbox{as\ } x_3 \rightarrow -\infty.
\end{array}
\right. 
\end{equation}

\begin{figure}
    \centering
\includegraphics[width = 0.7\linewidth]{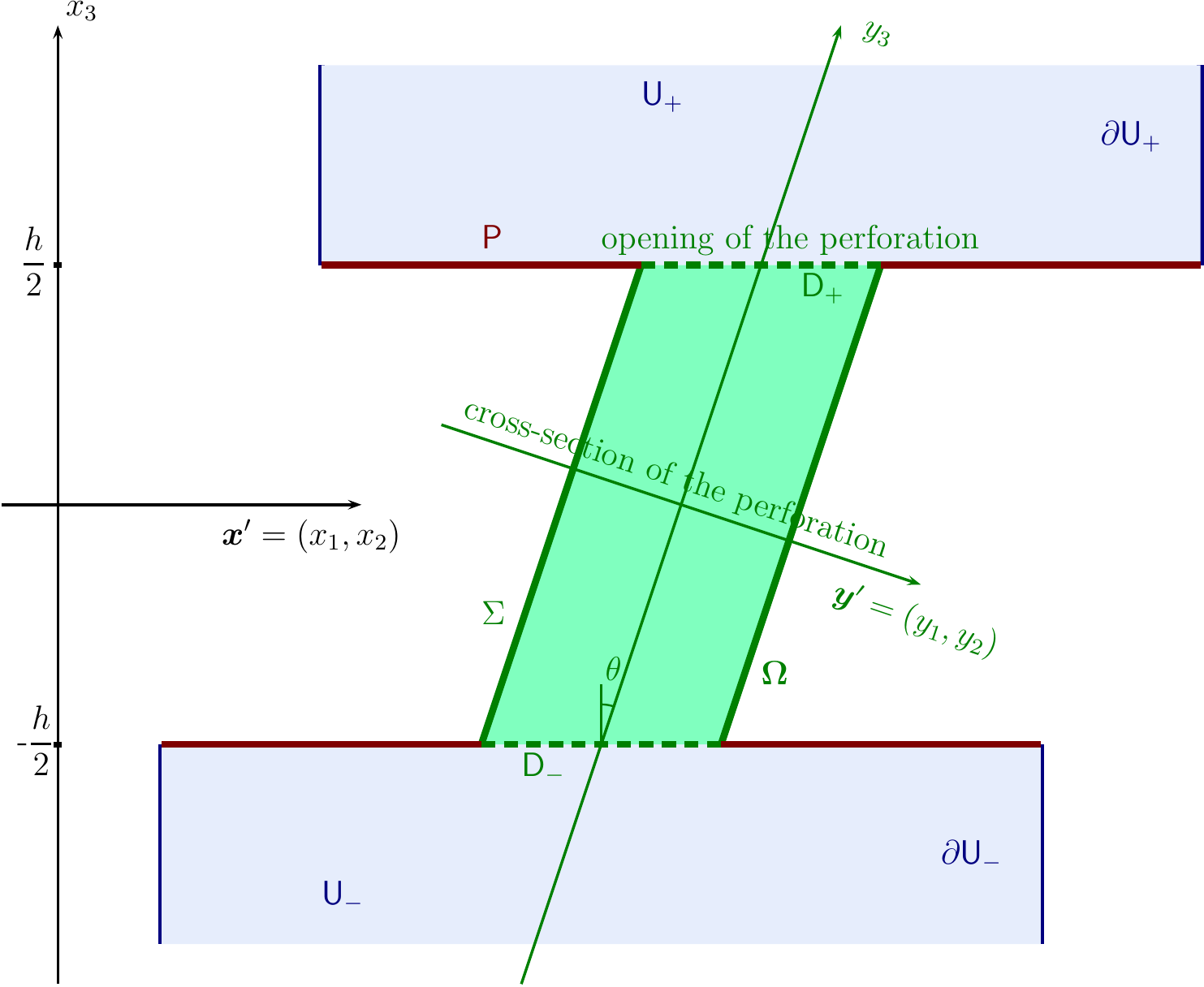}
    \caption{Schematic view of the lattice unit cell. The perforation is shown in green, the plate in red and the fictitious boundaries in blue. }
 \label{fig:domain}
  \end{figure}

\begin{figure}
    \centering
\subfloat[Opening (elliptical)]{\includegraphics[width = 0.3\linewidth]{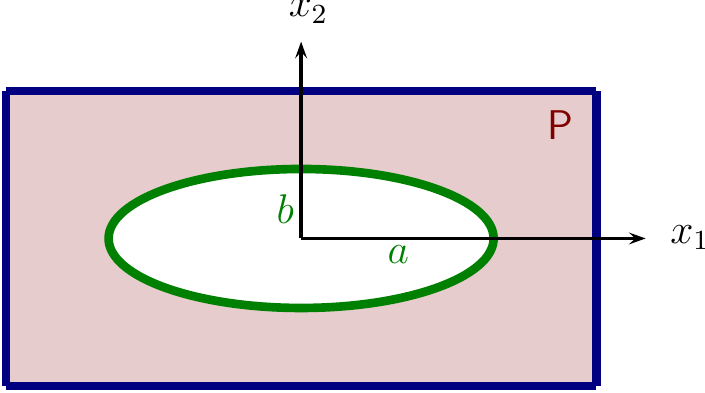}  \label{fig:opening}
}
\quad
  \subfloat[Elliptical cross-section]{\includegraphics[width = 0.33\linewidth]{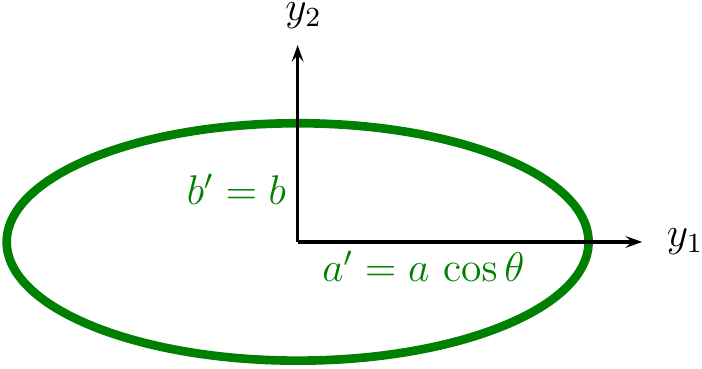} \label{fig:sectionell}
}
\quad
\subfloat[Circular cross-section]{ \includegraphics[width = 0.273\linewidth]{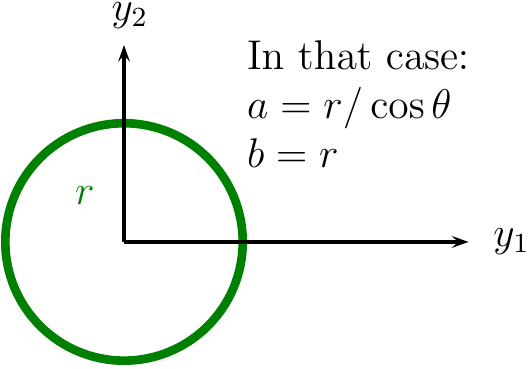} \label{fig:sectioncirc}}
    \caption{Opening and cross-section of the perforation.}
 \label{fig:sections}
  \end{figure}

The perforated plate contains a two-dimensional doubly-periodic lattice of apertures
\begin{equation}\label{eq:latticeproperties}
  \bxi_m=m_{1}\bxi_{1}+m_{2}\bxi_{2},
\end{equation}
where $m=(m_1,m_2)$ is a pairing of two integers and $\bxi_{i}$ are the periodicity vectors of the lattice (see figure \ref{fig:plaque-schema}). The low porosity assumption reduces the characterization of the acoustic behavior of each perforation as if it were isolated. As a result, the periodicity of the plate will not be used until section \ref{section:reseau}. The surface of the unit lattice cell is $A=\norm{\bxi_1\times \bxi_2}$ and the spacing between two successive perforations is denoted by $L=\max \left(|\bxi_{1}|,|\bxi_{2}| \right)$. All apertures have the same shape: they can be tilted, or untilted, their openings are elliptical while their cross-section can be elliptical or circular (see figure \ref{fig:domain} and  figure \ref{fig:sections}). This configuration is used for most industrial problems. In cooling liners used in combustion chambers, perforations are not perpendicular to the plate, but are tilted at an angle to the $x_3$ axis that is often around $60$\textdegree $ ~$ \citep[see][]{2007_AIAA_Eldredge, 2011_ASME_Andreini}. The general case of a tilted perforation with an elliptical opening and an elliptical cross-section (see figures \ref{fig:opening} and \ref{fig:sectionell}) is addressed in this paper. This case is representative of actual perforations. Especially, when a drill bit of radius $r$ is used for boring holes tilted by $\theta$ degrees from the vertical, it results in a perforation which cuts any parallel plane to the plate along an ellipse of major axis (in the tilt direction) $2a=2r/\cos\theta$ and minor axis (in the transversal direction) $2b=2r$ (see figures \ref{fig:opening} and \ref{fig:sectioncirc}).

 The acoustic wavelength $\lambda = 2\pi/\kappa$ is assumed to be large in comparison with the aperture size, however the spacing parameters $|\bxi_{1}|$ and $|\bxi_{2}|$ are only required to be less than half 
 a wavelength:  \citep[see][]{bendali2012}
\begin{equation}
  \label{eq:2.6}
  L < \lambda/2.
\end{equation}

 This ensures that only the fundamental mode, corresponding to $m_1=m_2=0$, is propagating for any incident wave. All the short-range interaction effects between the perforations can be neglected. For shorter wavelengths such that $L> \lambda/2$, the plate can no longer be assumed to react as an homogenized surface, whose reflecting properties are characterized by an effective impedance.

As shown in figure \ref{fig:domain}, each aperture of the perforated plate consists of a bounded domain  $\Omega$ located between the two planes $x_3 = -h/2$ and $x_3=h/2$. The tilting angle of the perforation is denoted by $\theta$.

The lower and upper openings of the perforation are denoted by $\sD_-$ and $\sD_+$ and its lateral part by $\Sigma$. It is convenient to consider the decomposition of domain $\sU$ within the unit cell into non-overlapping domains, filled by the fluid  consisting of, respectively, the lower and upper semi-infinite tubes $\sU_-$ and $\sU_+$ and the hole $\Omega$. The lateral parts of the boundaries of, respectively, $\sU_-$, $\sU_+$ and $\Omega$ will be denoted as $\p \sU_-$, $\p \sU_+$ and  $\p \Omega$. Finally, in Figure \ref{fig:sections}, $\sP$ stands for the part of the upper and lower surfaces of the plate within the unit cell. One can thus notice that the horizontal boundaries of  $\sU^+$  and  $\sU^-$ consist of the corresponding part of $\sP$ and $\sD^+$ and $\sD^-$ respectively.

Because the acoustic wavelength is large as compared with the aperture size, each aperture can be considered as acoustically compact. This means that the local motion through an aperture can be assumed incompressible: 
\begin{equation}
\bnabla \bcdot \bv=0 \quad \mbox{in\ }  \Omega,
\end{equation}
where $\bv = (v_1,v_2,v_3)$ is the acoustic velocity field. As a result, the governing equations (\ref{pb:Helmholtz})  can be expressed in terms of a potential:
 \begin{equation} \label{pb:potential}
\left\{
\begin{array}{ll}
\displaystyle
  \nabla^2 p = 0 &\quad \mbox{in\ }  \sU\\
   \p_\bn p = 0 &\quad \mbox{on the rigid walls} \\
 \lim_{\norm{\bx} \rightarrow \pm\infty} p = P^{\pm},
\end{array}
\right. 
\end{equation}
where $P^+$ and $P^-$ are limiting values of pressure at big distance from the aperture. By identification with (\ref{eq:pi}), (\ref{eq:pr}) and (\ref{eq:pt}), we have $P^+ = p_i+p_r$ and $P^- = p_t$. Under this assumption, the Rayleigh conductivity of the aperture is defined as \citep{Howe_book}:
\begin{equation}
K_R=\frac{\ci \omega \rho_0 Q}{P^+-P^-},
\label{eq:Ray-cond-Q}
\end{equation}
where $\rho_0$ denotes the mean density of the fluid and $Q$ is the volume flux through the aperture in the $x_3$ direction : $Q= \int_{D_+} v_3 \: \dif x_1 \dif x_2=\int_{D_-} v_3\: \dif x_1 \dif x_2$. Using the linearized momentum equation, the Rayleigh conductivity can be expressed as a function of pressure only:
\begin{equation}
K_R=\frac{1}{P^+-P^-}\int_{D_+} \p_{x_3} p \: \dif \bx'.
\label{eq:Ray-cond-p} 
\end{equation}

\section{Rayleigh conductivities of tilted elliptical perforations}\label{section:cond}

In this section, we propose a rigorous method for estimating the Rayleigh conductivity of thick plates with elliptical tilted perforations.
As in the methods of  \citet{rayleigh} and \citet{Howe_book}, the main ingredient is a variational theory linked to minimization of kinetic energy, but
 here both the dual Dirichlet and Kelvin variational principles come into play. This theory is found in Lagrangian and Hamiltonian mechanics and is of great importance in several fields, especially in mathematical modelling \citep[Chap 4, Sect 9]{couranthilbert} to obtain lower and upper bounds for potential and kinetic energies, in optimization \citep{luenberger1997optimization} to derive the dual formulations, and in numerical analysis  \citep[Chap 1]{ciarlet}  for assessing the validity or the accuracy of a numerical solution. 

 The estimates obtained by this method match those of \citet{Howe_book} for cylindrical and conical untilted perforations, but it can also be extended for other geometries.
 \citet{laurens2012} applied it to thick plates with perforations of circular openings, with or without a tilting angle, and to untilted conical perforations. 
 In this paper, we focus on tilted perforations with elliptical openings, which were not covered in \citet{laurens2012} but are more relevant since
  they are representative of actual perforations obtained by conventional drilling. Since the openings are no longer circular in shape, the computation of 
  the integral of the source density of an elliptical disk (which is mandatory for applying Kelvin and Dirichlet principles and obtaining estimates of the  Rayleigh conductivity) is a new issue addressed here. In \ref{sec:brieftheory}, we summarize some key points of \citet{laurens2012}. 
 Then, in \ref{sec:tiltedelliptical}, we focus on the case of tilted elliptical configurations, and describe a new technique for establishing bounds of the Rayleigh conductivity and of end corrections.
 Finally, in \ref{sec:derivativeconf}, it is shown that the estimates obtained for simpler geometries \citep[in, for example,][]{Howe_book, laurens2012}, can easily be considered as specific instances of the general case dealt with here.

\subsection{A brief description of the theory}\label{sec:brieftheory}

In the following section, the computations are done for elliptical openings, and require new and more complex tools than in \cite{laurens2012}. Technical details of the calculus  are provided in appendix \ref{sec:calculs}.
First of all, the boundary-value problem (\ref{pb:potential}) is converted to a more condensed form through the change of variable $p'=\left(p-(P^+ +P^-)/2\right)/\left(P^+-P^-\right)$. For simplicity, in what follows the prime will be removed. This yields:
 \begin{equation} \label{pb:potential-adim}
\left\{
\begin{array}{ll}
\displaystyle
  \nabla^2 p = 0 & \mbox{in\ }  \sU\\[0.125em]
   \p_\bn p = 0 &\mbox{on the walls} \\[0.125em]
 p^\pm =\lim_{x_3 \rightarrow \pm \infty} \: p = \pm \: 1/2 &\mbox{on\ }  \sU_\pm.\\
\end{array}
\right. 
\end{equation}

The Rayleigh conductivity (\ref{eq:Ray-cond-p}) can now be defined as:
\begin{equation}\label{eq:Ray-cond-p-adim} 
K_R=\int_{D_+} \p_{x_3} p \: \dif \bx'.
\end{equation}

Consider the functional 
\begin{equation}
  \sJ(\psi,\bq) =  \int_{U}\left| \bq(\bx)-\bnabla \psi(\bx)\right|^2 \: \dif \bx,
\end{equation}
where $(\psi,\bq)$ run over the spaces of functions $\sH_{1/2}(\sU)$ and $\sW(\sU)$,  $\sH_{1/2}(\sU)$ being the space of functions $\psi$ so that $\bnabla \psi$ is square integrable in $\sU$ with $\lim_{x_3\rightarrow  \pm \infty}\psi =\pm 1/2$ and  $\sW(\sU)$  the space of square integrable vectorial functions $\bq$ so that $\bnabla\bcdot\bq=0$ in $\sU$ and $\bq \bcdot\bn=0$ on the boundary of $\sU$. The functional $\sJ(\psi,\bq)$ can also be written in the form:
\begin{equation}
  \sJ(\psi,\bq) =  \int_{\sU}\left| \bq(\bx)\right|^2 \: \dif \bx + \int_{\sU}\left| \bnabla \psi(\bx)\right|^2 \: \dif \bx - 2\int_{\sU} \bnabla\psi\bcdot\bq(\bx) \: \dif \bx.
\end{equation}

\citet{laurens2012} proved that  the following relation holds, thanks to the properties of the spaces of functions $\sH_{1/2}(\sU)$ and $\sW(\sU)$:
\begin{equation}
 \int_{\sU} \bnabla\psi\bcdot\bq(\bx) \: \dif \bx=\int_{D^+} q_3 \: \dif \bx'=\int_{D^-} q_3  \: \dif \bx'.
\label{eq:func-equality}
\end{equation}

Thus, $\sJ$ can be decomposed into $ \sJ(\psi,\bq)\;=\;\sJ_1(\psi)-\sJ_2(\bq)$ with
\begin{equation}
\sJ_1(\psi)\;=\;\int_{\sU} \Big| \bnabla \psi(\bx)\Big|^2\: \dif \bx  \mbox{ and } \sJ_2(\bq)\;=\; 2\int_{D_+} q_3(\bx) \: \dif \bx'-\int_{\sU} |\bq(\bx)|^2\:\dif \bx.
\end{equation}

The Kelvin and Dirichlet principles are based on performing this decomposition and considering the dual expressions for the energy
 $\sJ_1(\psi)$ and $\sJ_2(\bq)$.

Obviously, as $p$ is the solution to system (\ref{pb:potential-adim}), it belongs to $\sH_{1/2}(\sU)$, $\bnabla p$ belongs to $\sW(\sU)$ and $\sJ(p,\bnabla p)=0$. 

Moreover, it is easy to see from equations (\ref{eq:Ray-cond-p-adim}) and (\ref{eq:func-equality}) that $\sJ_1(p)=K_R$, which directly leads to $K_R=\sJ_1(p)=\sJ_2(\bnabla p)$. It can then be established that $\sJ_1$ reaches its unique minimum at $p$ and $\sJ_2$ its unique maximum at $\bnabla p$, which leads to the following lemma \citep[see][Prop. (3.1)]{laurens2012}.

\begin{lemma}\label{prop:3.1} The same Rayleigh conductivity $K_R$ can be obtained from one of the following expressions:
 \begin{systeme}\label{syst:22}
 \textbf{Dirichlet principle: }&
 K_R\;=\;\min_{\psi\in\sH_{1/2}}\: \sJ_{1}(\psi)\\
& \mbox{ with }\sJ_1(\psi)\;=\; \int_{\sU} \Big| \bnabla \psi(\bx)\Big|^2\: \dif \bx\\
 \textbf{Kelvin principle: }&
 K_R\;=\;\max_{\bq\in\sW} \: \sJ_{2}(\bq)\\
&\mbox{ with }\sJ_2(\bq)\;=\; 2  \int_{D^+} q_3(\bx) \; \dif \bx'-\int_{\sU} |\bq(\bx)|^2\:\dif \bx.
 \end{systeme}
\end{lemma}
This lemma shows that lower and upper bounds for Rayleigh conductivity can be obtained for any arbitrary geometry of the perforation, 
provided that  the functionals $\sJ_1(\psi)$ and $\sJ_2(\bq)$ can be  explicitly evaluated. 

In \citet{laurens2012}, the Rayleigh conductivity was calculated for untilted cylindrical and conical perforations, as well as for a tilted perforation with circular openings.
 However, in an industrial context, the opening of the tilted perforation is elliptical since the perforation is bored with a cylindrical drill bit. 
In this paper, we thus consider the most typical case of tilted perforations with elliptical openings and give results both for the Rayleigh conductivity bounds and the end corrections, which are the commonly used parameters in the acoustics and fluid mechanics context.

\subsection{Application to tilted elliptical configurations}\label{sec:tiltedelliptical}

 Bounds for the Rayleigh conductivity can be derived from Lemma \ref{prop:3.1} by a suitable use of test functions $\psi$ and $\bq$  minimizing
  (respectively maximizing) the functional $\sJ_1$ (respectively $\sJ_2$) of (\ref{syst:22}). To do so, auxiliary problems arising from
  the elliptical tilted aperture have to be solved. 

The perforation $\Omega$ with a $0 \leq \theta <\pi/2$  tilting angle (see figure \ref{fig:domain}) 
\begin{equation}\label{eq:3.19}
\Omega=\Big\{\bx\in\mathbb{R}^3\;| \; (x_1, x_2)\in\sA \mbox{ and }  -h/2 < x_3 < h/2 \Big\},
\end{equation}
intersects any plane parallel to the plate of thickness  $h$  along an ellipse $\sA$ (see figure \ref{fig:opening}) of semi-major axis $a$ and semi-minor axis $b$, described for a fixed value of $x_3$ by
\begin{equation}
 \sA=\left\{\bx' \in\mathbb{R}^2\;\left|  \;  \frac{x_1^{2}}{a^{2}}+\frac{x_2^{2}}{b^{2}} < 1  \right. \mbox{ with  } a>b  \right\}.
 \end{equation}

 The eccentricity of the ellipse is denoted by  $\epss = \sqrt{1 - b^2/a^2}$ in which are expressed the complete elliptic integral of the first kind $K(\epss)$ and second kind $E(\epss)$:
\begin{equation}
  \label{eq:Keps}
  K(\epss) = \int_0^{\pi/2} \left( 1-\epss^2 \sin^2 \phi \right)^{-1/2} \: \dif \phi \: \mbox{ and } \:
   E(\epss) = \int_0^{\pi/2} \left( 1-\epss^2 \sin^2 \phi \right)^{1/2} \: \dif \phi.
 \end{equation}

Noticing that $ K(0)= \pi/2$, we also consider the integral $D(\epss)$
 \begin{equation}
   \label{eq:Deps}
   D(\epss)  = \int_0^{\pi/2} \frac{ \sin^2\phi  }{\sqrt{ 1-\epss^2 \sin^2 \phi } } \: \dif \phi =  \frac{ K(\epss) - E(\epss)}{\epss^2},
 \end{equation}
satisfying $D(0)= \pi/4$. To find the suitable test functions $\psi$ and $\bq$, we compute the potential $f_{\epss}$ generated by a source density function
 $\rho_{f_{\epss}}$ defined on the elliptical hole $\sA$ :
\begin{equation}
  \label{eq:electrostaticpotential}
  f_{\epss}(\bx',x_3) = \frac{1}{2 \pi} \int_\sA \frac{ \rho_{f_{\epss}}
    (\by)}{|\bx' - \by|} \dif \by, \quad  \bx' \in \sA.
\end{equation}

 Conveniently, we consider the following local coordinate systems of respectively $\sU^\pm$:
 \begin{equation}
  \bx_\pm=(\bx'_\pm,x_3^\pm) \quad \mbox{ with }  x_3^+ \geq 0 \mbox{ and } x_3^- \leq 0.
 \end{equation}

This potential is the solution of:
\begin{equation}\label{eq:pbpotential}
\left\{  \begin{array}{ll}
 \nabla^2 f_{\epss}(\bx_\pm)=0,&\mbox{ for }x_3^\pm \neq 0,\\[0.5em]
 \partial_{x_3^\pm}f_{\epss}(\bx'_\pm,0)=0,&\mbox{ for } \bx'_\pm \not \in\sA .
  \end{array} \right.
\end{equation}

Depending on the choice of boundary condition on the elliptical hole $\sA$, \ie the choice of the source density function, various expressions for the potential can be obtained.

First, we consider the case for which the potential verifies a Dirichlet condition on the elliptical hole. As  the aperture is tilted with an angle $\theta$,
 the projection of the potential on $\sA$ is a polynomial of degree $1$. 
We separate it into two components, $w_{\epss}$ and $t_{\epss}$, solutions of (\ref{eq:pbpotential}) with the following boundary conditions:
\begin{equation}\label{eq:30}
 w_{\epss}(\bx',0)= 1 \quad \mbox{ and } \quad  t_{\epss}(\bx',0)=x_1 \quad \mbox{ for } \bhx \in\sA.
\end{equation}

 For $\alpha \in\mathbb{R}$,  the test function $\psi$ to use in Lemma \ref{prop:3.1}  is defined as
\begin{equation}\label{eq:psi}
\psi(\bx) = 
\left\{  \begin{array}{ll}
\pm \frac{1}{2} \mp \left(1 -  2\alpha\right) \: w_{\epss}(\bx_\pm)
+ \left(\mu_\alpha\: \sin \theta \:/\:h \right) \: t_{\epss}(\bx_\pm),&\mbox{ for } \pm x_3^\pm >0,\\
 -\alpha+  \left( \mu_\alpha \:/\:h \right)\Big( (x_3 + h/2 ) \:\cos \theta+x_1\:\sin \theta \Big), &\mbox{ for }  x_3^+ x_3^-  < 0 ,
\end{array} \right.
\end{equation}
with $\mu_\alpha= 2 \alpha \: \cos \theta$.
The functional $\sJ_1$ involved in the Dirichlet principle of Lemma \ref{prop:3.1} is thus 
\begin{equation}  \label{eq:j1}
  \begin{array}{ll}
\ds \sJ_1(\psi) = & \ds 2  ( 1- 2 \alpha)^2 \int_{\small \pm x_3^\pm>0}  \big| \bnabla w_{\epss}(\bx) \big|^2\;\dif \bx \\
   &  + \ds \: 2 \:  \frac{\mu^2_\alpha\: \sin^2 \theta}{h^2}\int_{\small\pm x_3^\pm>0}  \big| \bnabla t_{\epss}(\bx) \big|^2 \;\dif \bx 
 + \frac{\mu^2_\alpha}{h^2} ( \pi a b \: h).
  \end{array}
\end{equation}

By using the chosen expressions (\ref{eq:integrals}) of the integrals (see appendix \ref{sec:calculs}), we then get
\begin{equation}
  \begin{array}{ll}
 \ds  \sJ_1(\psi)=& \ds 4 a  \frac{ K(0)}{K(\epss)} -  \ds \left(  16 a \frac{K(0)}{K(\epss)}  \right) \alpha\\
 &\ds  +   \left( 16 a  \frac{K(0)}{K(\epss)} + \frac{64a^3}{3 h^2}\cos^2 \theta \sin^2 \theta   \frac{D(0)}{D(\epss)} +   \frac{4 \pi a b}{h}   \cos^2 \theta \right) \alpha^2.
  \end{array}
\end{equation}
The minimal value for this function $\sJ_1(\psi)$ is reached for
\begin{equation}
  \label{eq:alphaj1}
   \alpha =  8 a \frac{ K(0)}{K(\epss)}  \left/   \left( 16 a   \frac{K(0)}{K(\epss)} + \frac{64 a^3}{3 h^2}  \cos^2 \theta \sin^2 \theta \frac{D(0)}{D(\epss)} +   \frac{4 \pi a b }{h}\cos^2 \theta \right) \right.,
\end{equation}
which leads to the upper bound for the Rayleigh conductivity.

For finding the lower bound, we have to consider the potential $z_{\epss}$ solution of (\ref{eq:pbpotential}) with a Neumann boundary condition on the elliptical hole:
\begin{equation}\label{eq:25}
\partial_{x_3^\pm}z_{\epss}(\bx', x_3^\pm = 0)=\pm 1/2 \mbox{ for } \bhx \in\sA .
\end{equation}

The suitable test function for $\sJ_2$ is
\begin{equation}\label{eq:31}
\bq(\bx) = 
\left\{  \begin{array}{ll}
 \pm (\beta \: / \:\pi ab ) \: \bnabla z_{\epss}(\bx_\pm)&\mbox{ for } \pm x_3^\pm> 0,\\[0.5em]
 (\beta \: /\: \pi ab ) \: \Big(\be_3+\tan \theta \: \be_1\Big)&\mbox{ for } x_3^+ x_3^- < 0,
\end{array} \right.
\end{equation}
depending on a real constant $\beta$.

The functional involved in the Kelvin principle of Lemma \ref{prop:3.1} is then:
\begin{equation}
  \label{eq:j1}
  \sJ_2(\bq)= 2 \beta -  \frac{\beta^2}{\pi^2\;a^2 b^2}  \left(  2 \int_{\small \pm x_3^\pm>0}  \big| \bnabla z_{\epss}(\bx) \big|^2\;\dif \bx   + ( 1 + \tan^2 \theta) \pi a b \: h \right) .
\end{equation}

The analytical expression of the integral term is calculated in the appendix \ref{sec:calculs}, equation (\ref{eq:25bis}), which leads to:
\begin{equation}
  \label{eq:j1bis}
  \sJ_2(\bq)= 2 \beta -  \frac{\beta^2}{\pi^2\;a^2 b^2}  \left(  \frac{16}{3} ab^2  \frac{K(\epss)}{K(0)} + ( 1 + \tan^2 \theta) \pi a b \: h \right) .
\end{equation}

The maximal value of this function is reached for
\begin{equation}
  \label{eq:beta}
  \beta = \pi\;a b \left/  \left(  \frac{16\, b}{3 \pi } \frac{K(\epss)}{K(0)}  + ( 1 + \tan^2 \theta)  h \right) \right..
\end{equation}

Finally, the bounds of the Rayleigh conductivity are 
 \begin{equation}
  \label{eq:ellipseinc}
   \frac{\pi ab}{  \frac{h}{\cos^2 \theta} + \frac{16b}{3\pi}\frac{K(\epss)}{K(0)} }\leq  K_R \leq \frac{\pi ab}{  \frac{\pi b}{2}\frac{K(\epss)}{K(0)} + \frac{h}{\cos^2 \theta} \left/ \left(1+\frac{16a^2}{3 \pi b h}\frac{D(0)}{D(\epss)} \sin^2 \theta\right) \right. }.
\end{equation}

It is instructive to express these bounds in terms of the effective length of the perforation, meaning the length $l = s/K_R$, where $s$ is the area of the aperture. 
This length can be interpreted as the effective length of a `piston' of fluid of cross section $s$ involved in the motion through the aperture. 

For example, for a cylindrical untilted aperture of radius $r$, the conductivity of the channel alone is $\pi r^2/ h$, for which $l$ would be $h$. Thus, the Rayleigh conductivity can be expressed in terms of an end correction $l'$, defined by
\begin{equation}
  \label{eq:endcorrection}
  l =h+l',
\end{equation}
so that $l'$ is the amount by which $h$ has to be increased to account for the contributions of the openings. 

It should be noted that this end correction can  be directly linked to the specific acoustic impedance of the perforation, defined in Sect. \ref{section:notations} as:
\begin{equation}
 z=\frac{1}{\rho_0 c_0}\frac{P^+-P^-}{- v_3}.
\end{equation}

As proposed by \citet{2000_AIAAJ_Jing}, the impedance can indeed be decomposed into the impedance induced by the end correction $z_e=-\ci \kappa l'$ and the added thickness term $-\ci \kappa h$:
\begin{equation}
 z=z_e-\ci \kappa h=-\ci \kappa (l'+h).
\end{equation} 

 However, for tilted perforations, the end correction analogy is not as straightforward as for untilted configurations. Indeed, the end correction is related to the choice of the reference area $s$ so that $z = - \ci\kappa s /K_R$. Because of the definition of the volume flux involved in the Rayleigh conductivity (see (\ref{eq:Ray-cond-Q})), $s$ should be the surface area of the opening. But, in the literature \citep{2011_ASME_Andreini,2009_JCP_Mendez, 2007_AIAA_Eldredge},  the cross-section area is used instead, with a modified perforation length $h/\cos{\theta}$. Strictly speaking, it is a not justified choice. Equation (\ref{eq:ellipseinc}) allows a rigorous definition of the end correction for a tilted perforation.

Since the opening is elliptical,  the effective section is $s = \pi \, ab$ and the bounds for the effective length are
 \begin{equation}
   \label{eq:lellinc} \frac{\pi b}{2}\frac{K(\epss)}{K(0)} +  \frac{h}{\cos^2 \theta}\left/ \left(1+\frac{16a^2}{3 \pi b h}\frac{D(0)}{D(\epss)} \sin^2 \theta\right) \right.
\leq  l\leq \frac{16b}{3\pi}\frac{K(\epss)}{K(0)}+ \frac{h}{\cos^2 \theta}.
 \end{equation}

To define the end correction $l'$ similarly as in (\ref{eq:endcorrection}), we thus introduce the modified effective height $\widetilde{h}=h/\cos^2 \theta$, so that 
\begin{equation}
  \label{eq:endcorrectiontilted}
  l =\widetilde{h}+l'
\end{equation}
and
 \begin{equation}
   \label{eq:lellipseinc} \frac{\pi b}{2}\frac{K(\epss)}{K(0)} + g(h)
\leq  l'\leq \frac{16b}{3\pi}\frac{K(\epss)}{K(0)},
 \end{equation}
with $  g(h) =  -  \frac{16 a^2}{3 \pi b }\frac{D(0)}{D(\epss)}\tan^2 \theta \left/ \left( 1+\frac{16a^2}{3 \pi b h}\frac{D(0)}{D(\epss)} \sin^2 \theta\right) \right. \leq 0$.

It is worth noting that the effective height $\widetilde{h}$ is hence quite different from the intuitive effective thickness of the perforation $h/\cos{\theta}$ commonly used in the literature.

However, in the specific case of a circular cross-section of radius $r$ (see figure \ref{fig:sectioncirc}), the opening is elliptical with semi-major axis $a = r / \cos \theta$ and semi-minor axis $b=r$ (see figure \ref{fig:opening}). 
Thus, the eccentricity of the opening is $\epss=\sin \theta$. The opening area $s$ is $\pi a b=\pi r^2/\cos \theta$ and the bounds (\ref{eq:ellipseinc}) for the Rayleigh conductivity can be set as follows:
\begin{equation}  \label{eq:ellipseinc_bis}
\begin{array}{ll}
  \frac{\pi r^2}{ \frac{h}{\cos \theta} + \frac{16 r}{3\pi}\frac{K(\sin \theta)}{K(0)} \cos \theta} & \leq  K_{R} \\
&  \leq \frac{\pi r^2}{  \frac{\pi r}{2}\frac{K(\sin \theta)}{K(0)} \cos \theta  + \frac{h}{\cos 
\theta} \left/ \left(1+\frac{16 r}{3 \pi  h}\frac{D(0)}{D(\sin \theta)} \tan^2 \theta\right)\right.}
\end{array}\end{equation}

This means that if the cross-section $\pi r^2$ is considered as the reference section instead of the effective section $\pi a b=\pi r^2/\cos \theta$, 
the effective height, which then has to be considered in the definition of the end correction, is indeed the intuitive value $h/\cos \theta$ instead of $h/\cos^2 \theta$.

\subsection{Derivate geometrical configurations}\label{sec:derivativeconf}

In this subsection, we show that classical results for the Rayleigh conductivity and the end correction of simple geometries, as well as the results of \citet{laurens2012} for tilted perforations with circular openings are contained in  the general expressions (\ref{eq:ellipseinc}) and (\ref{eq:lellipseinc}). However, tapered holes \citep[see][chapter 4, for instance]{1999_NASA}, which may be obtained when perforated plates are manufactured by laser drilling, cannot be obtained from (\ref{eq:ellipseinc_bis}) valid only for cylindrical perforations characterized by identical openings on each side of the plate. However, these bounds for a conical untilted perforation with circular opening can be found in \citet{laurens2012}.

\subsubsection{Untilted perforation with an elliptical opening}

We set $\theta = 0$ in  (\ref{eq:ellipseinc}) to get:
\begin{equation}
  \label{eq:ellipsedroite}
  \frac{\pi ab}{ h + \frac{16b}{3\pi}\frac{K(\epss)}{K(0)} }\leq  K_{R} \leq \frac{\pi ab}{ h+\frac{\pi b}{2}\frac{K(\epss)}{K(0)}}.
\end{equation}

As explained by \citet{1969_JSV_Morfey}, the upper bound of equation (\ref{eq:ellipsedroite}) was found by Lord Rayleigh \citep[][chapter 1, section 306]{rayleigh} for an infinitely thin plate ($h=0$). The lower bound expression and the taking into account of the plate thickness is, to the best of the authors' knowledge, entirely new.

This result can once again be expressed in terms of an end correction. Using the effective section of the elliptical opening $s = \pi ab$, the bounds for this end correction are:
\begin{equation}\label{eq:lellipsedroite}
\frac{\pi b}{2}\frac{K(\epss)}{K(0)}\leq  l' \leq \frac{16b}{3\pi}\frac{K(\epss)}{K(0)}  .
\end{equation}

\subsubsection{Tilted perforation with a circular opening}

This case was considered by \citet[Thm $4.3$]{laurens2012}, but is not the most useful for actual applications, since the manufacturing of such perforations would require an elliptical drill bit. The Rayleigh conductivity can be deduced from (\ref{eq:ellipseinc}) in a straightforward way by taking  $\epss = 0$ (which avoids elliptic integrals):
 \begin{equation}
   \label{eq:cylinc}
   \frac{\pi r^2}{ \frac{h}{\cos^2 \theta} + \frac{16r}{3\pi} } \leq  K_R\leq \frac{\pi r^2}{
 \frac{h}{\cos^2 \theta} \left/ \left(1+\frac{16r}{3 \pi h} \sin^2 \theta\right)\right. +\frac{\pi r}{2}}.
 \end{equation}

Once again, this result shows that the tilting angle has a non-intuitive influence on the Rayleigh conductivity of the perforation, since 
the effective height to be considered in the end correction definition (\ref{eq:endcorrectiontilted}) is $\widetilde{h}=h/\cos^2 \theta$. The bounds for the end correction are given by:
 \begin{equation}
    \label{eq:lcylinc}
  f(h) + \frac{\pi r}{2}
\leq  l' \leq \frac{16r}{3\pi},
  \end{equation}
with $ f(h) = -  \frac{16r}{3\pi} \tan^2 \theta \left/ \left(1+\frac{16r}{3 \pi h} \sin^2 \theta\right) \right.$.

\subsubsection{Untilted perforation with a circular opening}\label{sec:untiltedcylinder}

The case of an untilted cylindrical perforation has already been addressed by \citet{Howe_book}. He established the following bounds for the Rayleigh conductivity:
\begin{equation} \label{eq:bounds-untiltedcylinder}
  \frac{\pi r^2}{h + 16 r \: / \:3 \pi} \leq K_R\leq \frac{ \pi r^2}{ h + \pi r / 2},
\end{equation}  
or, in terms of an end correction:
\begin{equation} \label{eq:lbounds-untiltedcylinder}
\frac{ \pi r}{ 2} \leq l' \leq  \frac{16 r}{3 \pi}
\end{equation} 
 
Identical values to those of Howe are obtained using (\ref{eq:ellipseinc}) with $\cos \theta = 1$ since the aperture is untilted and $K(\epss) = K(0)$ since the eccentricity of a disc is $\epss = 0$.

\subsubsection{Infinitely thin plate with elliptical or circular apertures}

This configuration has already been studied by \citet[formula $(5)$, section $306$]{rayleigh} and \citet[formula $(2.32)$, section $2$]{1973_JFM_Leppington}. They established $ K_R   ~= ~\pi a / K(\epss) $, which, in the limit case $\epss = 0$ ($a=r$ and $K(0)=\pi/2$) gives $K_R = 2r$; the same result was found by \cite{copson1947} by a direct computation of the flux integral on a circle. 

 Taking $h = 0$ in  (\ref{eq:ellipsedroite}) (untilted perforation with an elliptical opening, derived from (\ref{eq:ellipseinc}) by setting $\theta = 0$) then gives:
 \begin{equation}\label{eq:eq:bounds_nothickness}
   \frac{ 3 \pi^ 2 a }{ 32 K(\epss)}  \leq K_R\leq \frac{ \pi a }{ K(\epss)}.
 \end{equation}
and our upper bound  matches Rayleigh's formula exactly. For a circular aperture, doing the same with   (\ref{eq:bounds-untiltedcylinder}) leads to
 \begin{equation}\label{eq:eq:bounds_nothickness_circle}
   \frac{ 3 \pi^ 2 r }{ 16}  \leq K_R\leq 2r.
 \end{equation}

Once again, the upper bound is the analytical value of the Rayleigh conductivity of a circular perforation in an infinitely thin plate given by  \cite{rayleigh} and  \cite{copson1947}.

  \section{Numerical results}\label{sec:numericalresults}

In this section, the influence of the tilting angle and opening eccentricity on the Rayleigh conductivities and end corrections is illustrated for some configurations representative of 
perforated plates in combustors.

\subsection{Typical perforated plate}\label{sec:generic-perforated-plate}
The multiperforated plate under consideration is among those typically involved in combustion chambers \citep[see][for instance]{2011_ASME_Andreini, 2012_ASME_Andreini}. Such plates have quite a low porosity ($\approx 2\%$). Consequently the individual apertures can be assumed to be coupled only through long range interactions which means that the results of previous sections apply. The plate thickness is $h=2$mm and the characteristic size of the perforations is $r= 0.225$mm. This length will be  the semi-minor axis $b$ of the elliptical opening of the perforation (see figure \ref{fig:opening}). In the specific case of a hole bored with a drill bit of radius $r$, the cross-section of the tilted perforation is circular with the same radius, and the elliptical opening is characterized by $b=r$, $a=b/\cos\theta$ and $\epss=\sin\theta$ (see figure \ref{fig:sectioncirc}).
 
The reference values for the Rayleigh conductivity and the end correction will be respectively the upper and lower bounds of these parameters, corresponding to an untilted cylindrical hole (see section \ref{sec:untiltedcylinder}): 
\begin{equation}
 K_R^{\mathrm{ref}}=\frac{\pi r^2}{h+\pi r/2},  \quad l_{\mathrm{ref}}=\pi r/2
\end{equation}
In particular, we have $h=5.66 \, l_{\mathrm{ref}}$.  

\begin{figure}
  \centering  
 \includegraphics[width = 0.45\linewidth]{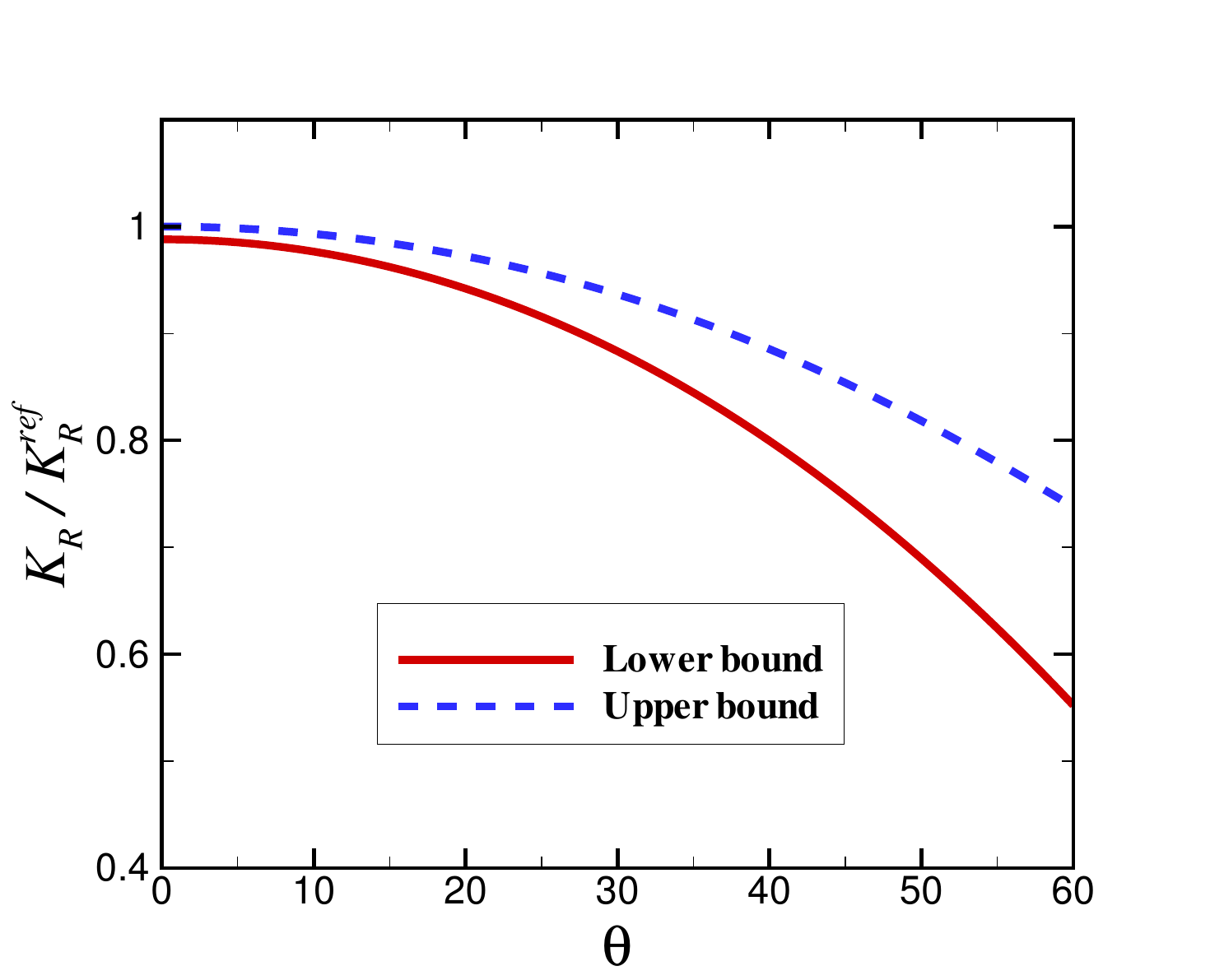}
 \caption{Scaled Rayleigh conductivities for the tilted perforation with an elliptical opening and a circular cross-section. The tilting angle $\theta$ is in degrees.}
  \label{fig:KR_incline}
 \end{figure}

First, the influence of the tilting angle is illustrated for the specific case of the tilted perforation with a circular cross-section. Bounds for the Rayleigh conductivity are given by
(\ref{eq:ellipseinc_bis}) and plotted in figure \ref{fig:KR_incline}, for a tilting angle varying from 0 to the conventional $60$\textdegree~value.

\begin{figure}
     \centering   
 
 \quad \subfloat[Scaled effective thickness]{\includegraphics[width = 0.45\linewidth]{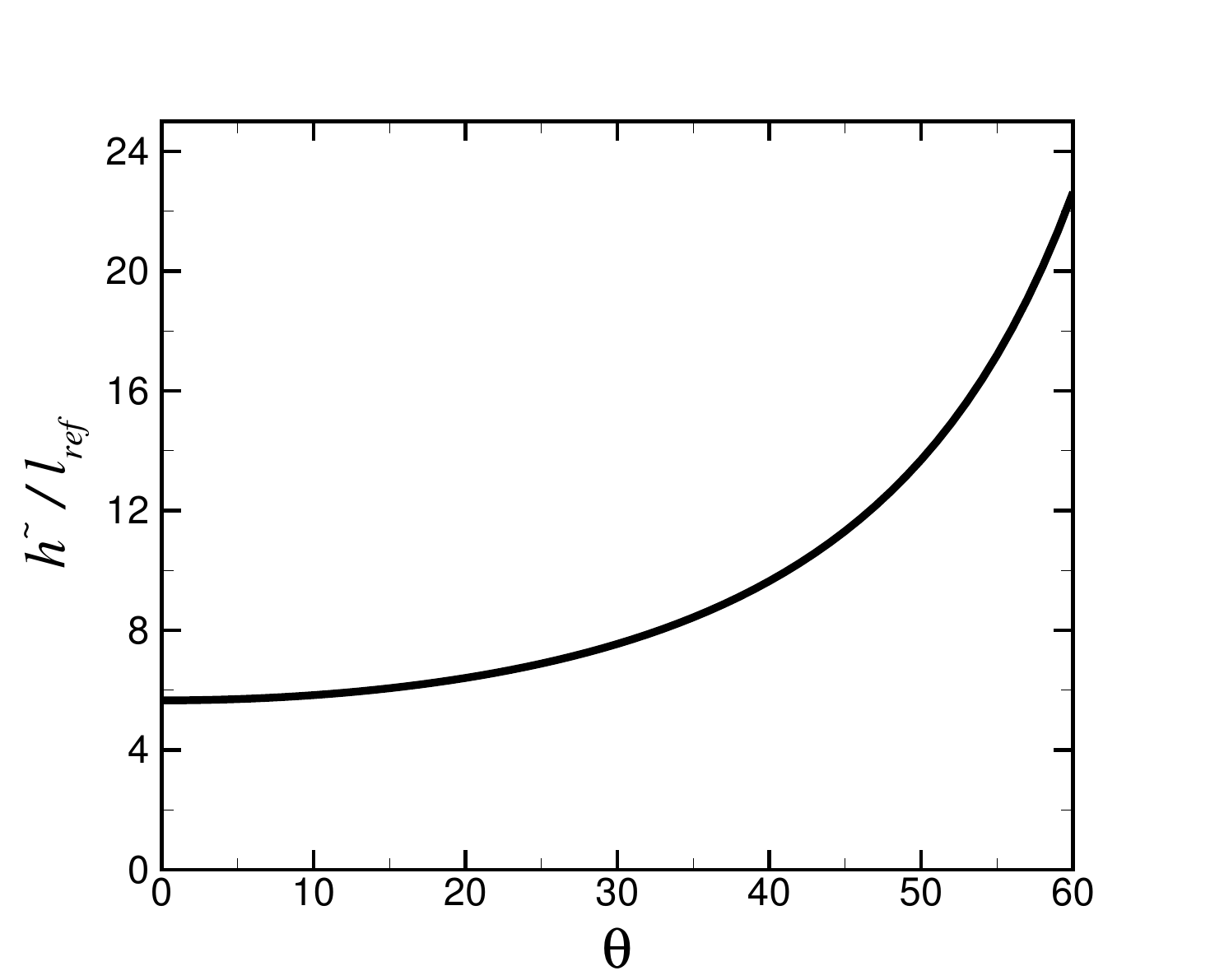} \label{fig:htilde}
 }
  \quad \subfloat[Scaled end correction]{\includegraphics[width = 0.45\linewidth]{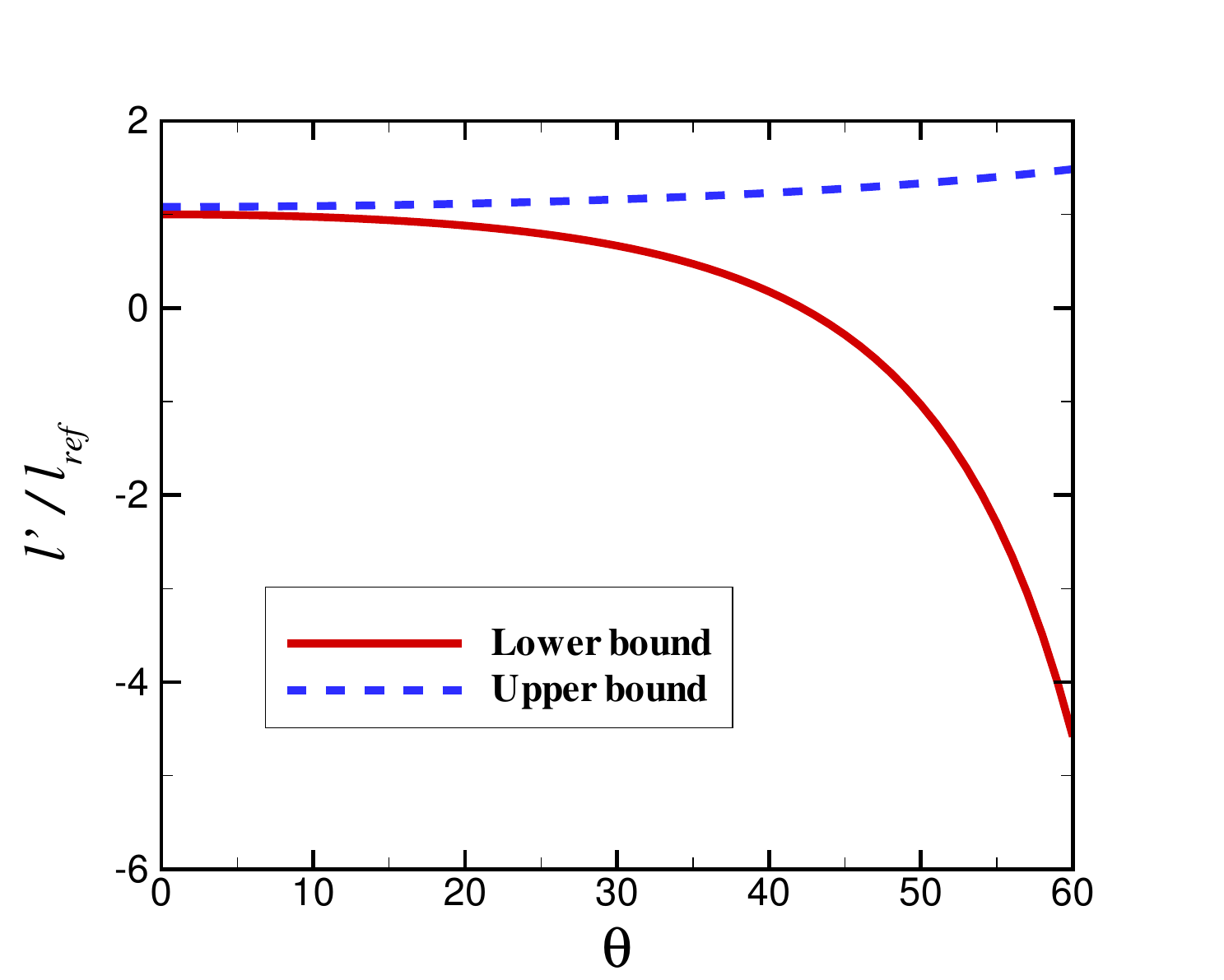} \label{fig:lprime_ellipse_inclinee}
 }
 \quad \subfloat[Scaled effective length]{\includegraphics[width = 0.45\linewidth]{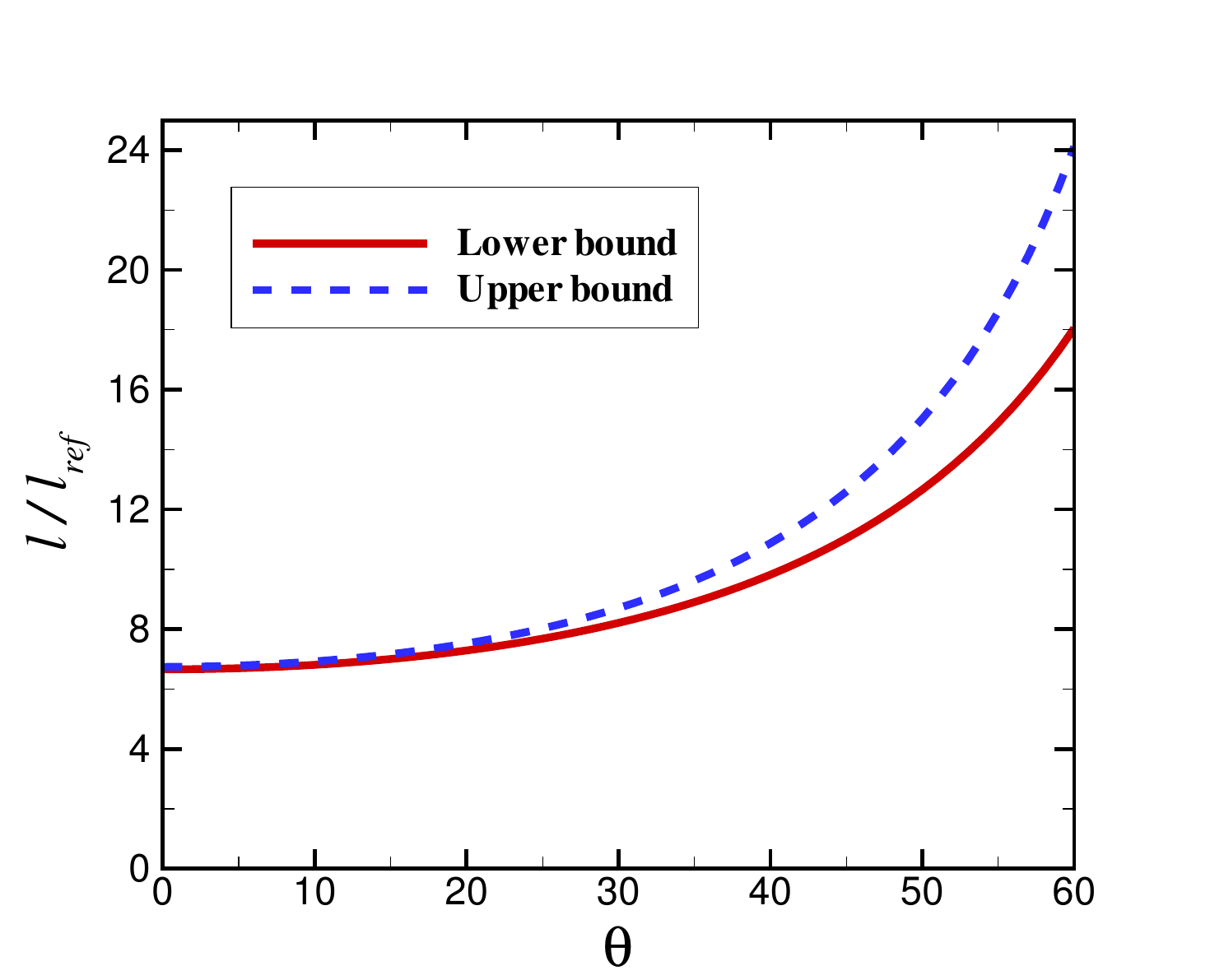} \label{fig:ltotal_ellipse_inclinee}
 }
    \caption{Scaled effective lengths and end correction for the tilted perforation with an elliptical opening and a circular cross-section. The tilting angle $\theta$ is in degrees. }
  \label{fig:l_incline}
 \end{figure}

A significant deviation from the untilted cylindrical case is observed at the largest tilting angles, which shows the importance of correctly taking the tilting effect into account. The effective length of the perforation is plotted in figure \ref{fig:l_incline}, as well as its decomposition into effective thickness and end correction (it should be remembered that $l=\widetilde{h}+l'$, see equation (\ref{eq:endcorrectiontilted})). A divergence between lower and upper bounds of the end correction is observed for large values of the tilting angle. However, this is not that significant for the effective length, because of the prominent role of the effective thickness $\widetilde{h}$ in this range of tilting angles.

\begin{figure}
  \centering  
 \includegraphics[width = 0.45\linewidth]{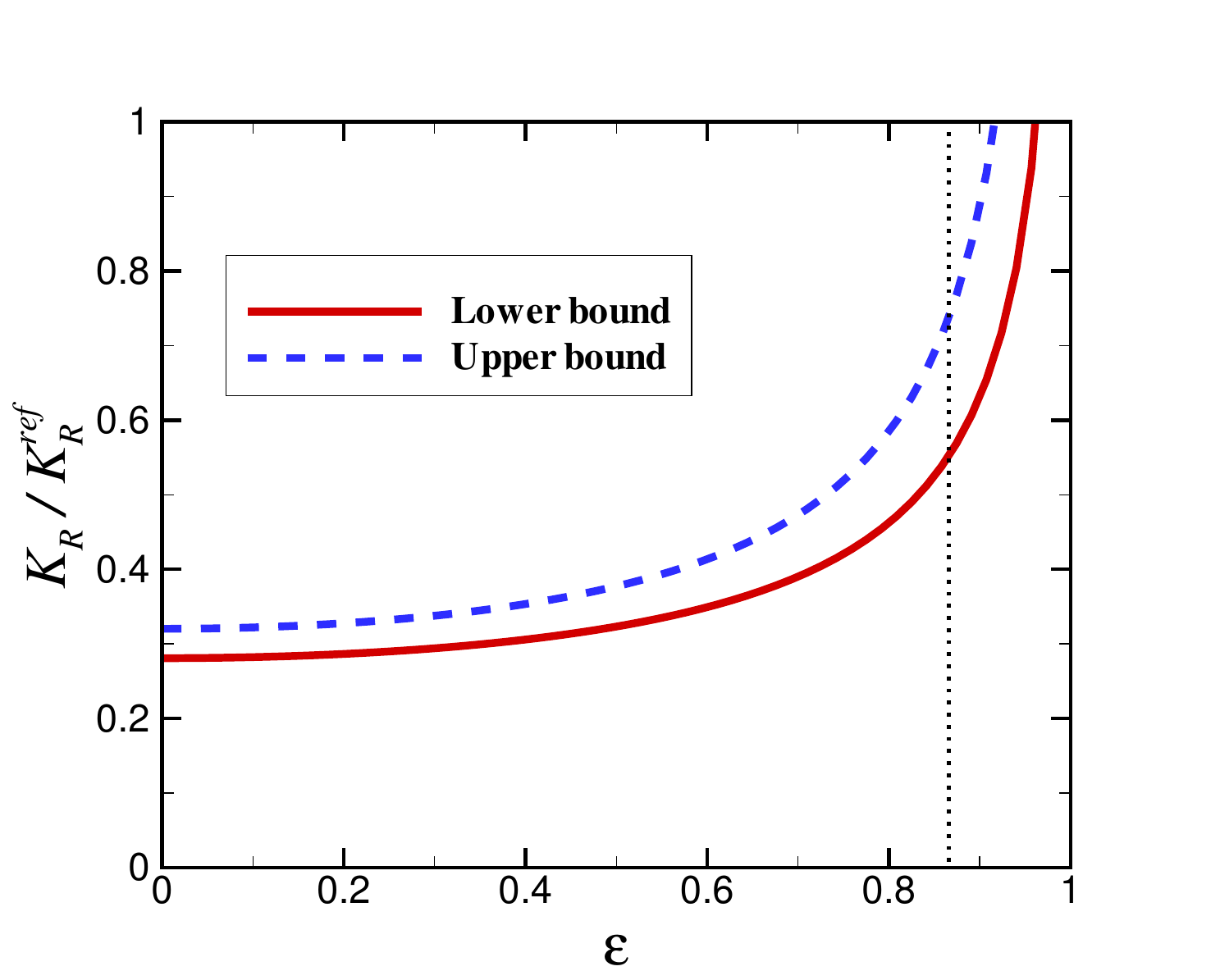}
 \caption{Scaled Rayleigh conductivities for an elliptical perforation with a $\theta=60$\textdegree~tilting angle. $\epss$ is the eccentricity of the elliptical opening. The dotted vertical line 
    corresponds to the eccentricity for which the cross-section of the perforation is circular.}
  \label{fig:KR_incline_epsilon}
\end{figure}

Then, the influence of the eccentricity of the elliptical opening is illustrated for a perforation with a $\theta=60$\textdegree~tilting angle, and a semi-minor axis  kept constant to $b=r$. Bounds for the scaled  Rayleigh conductivity given in (\ref{eq:ellipseinc}) are plotted in figure \ref{fig:KR_incline_epsilon}. The eccentricity corresponding to the specific case of a circular cross-section
(\ie  $\epss=\sin\theta=\sqrt{3}/2$) is denoted by the dotted vertical line.

\begin{figure}
     \centering   
 
 \quad \subfloat[Scaled effective thickness]{\includegraphics[width = 0.45\linewidth]{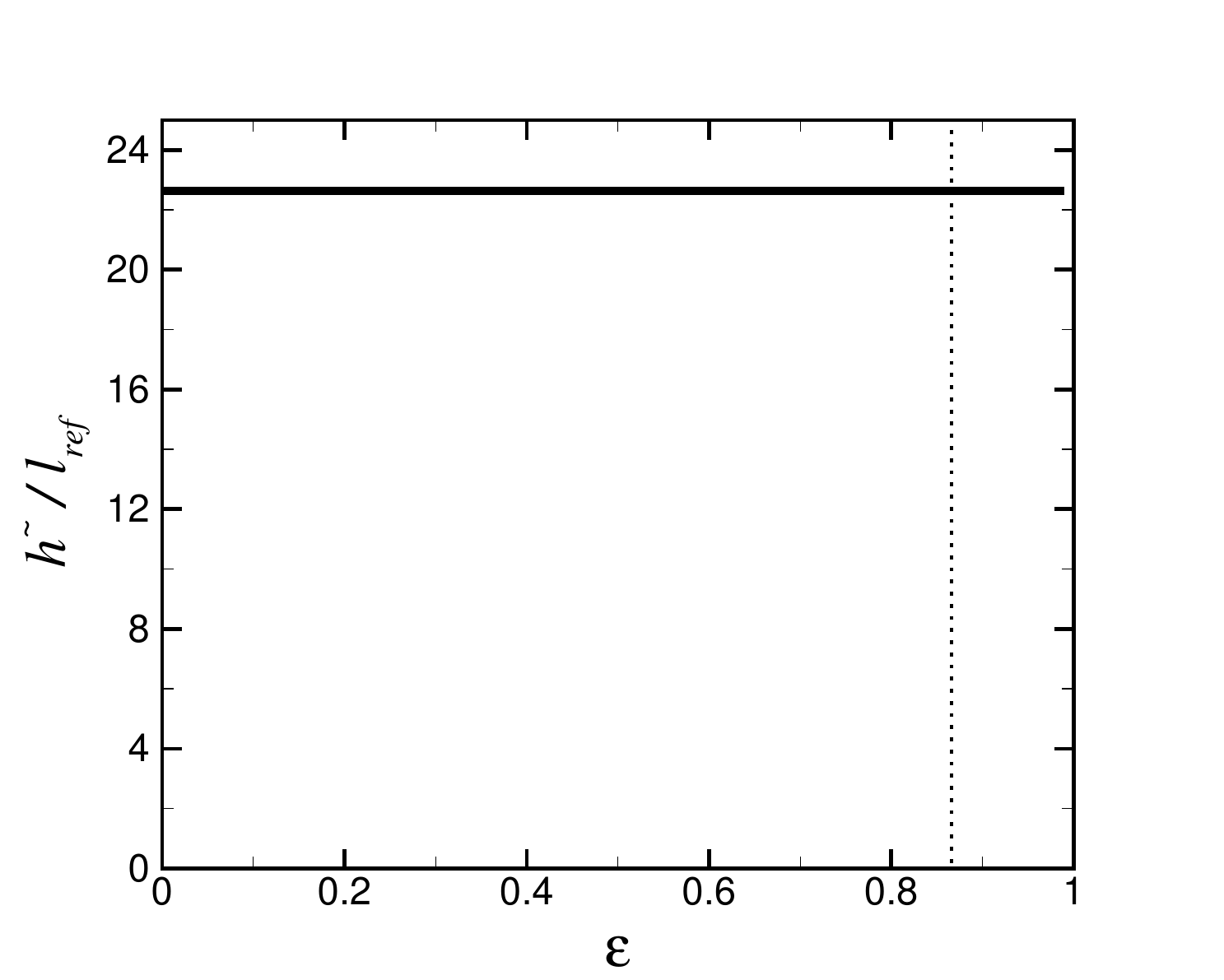} \label{fig:htilde_epsilon}
 }
  \quad \subfloat[Scaled end correction]{\includegraphics[width = 0.45\linewidth]{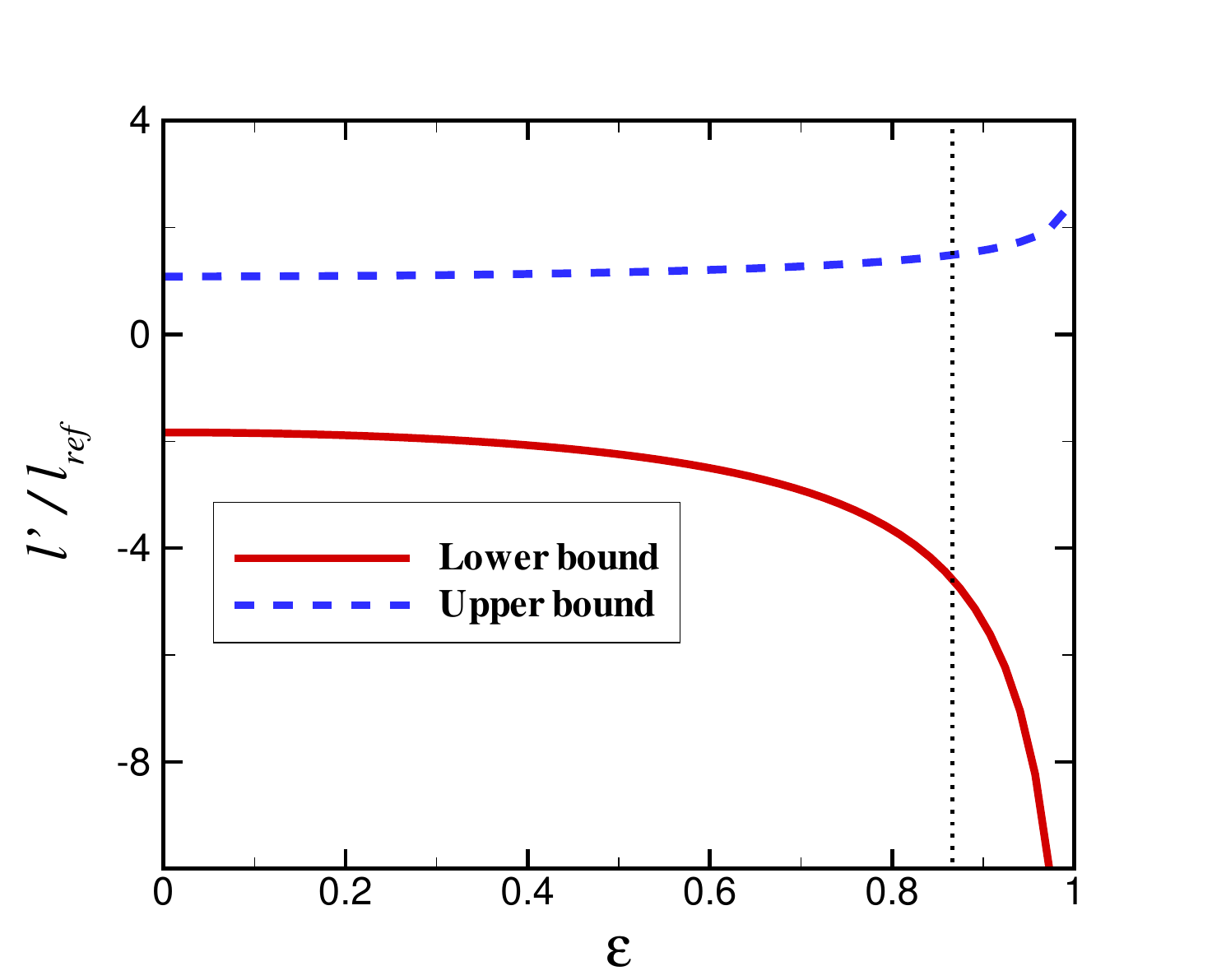} \label{fig:lprime_ellipse_inclinee_epsilon}
 }
 \quad \subfloat[Scaled effective length]{\includegraphics[width = 0.45\linewidth]{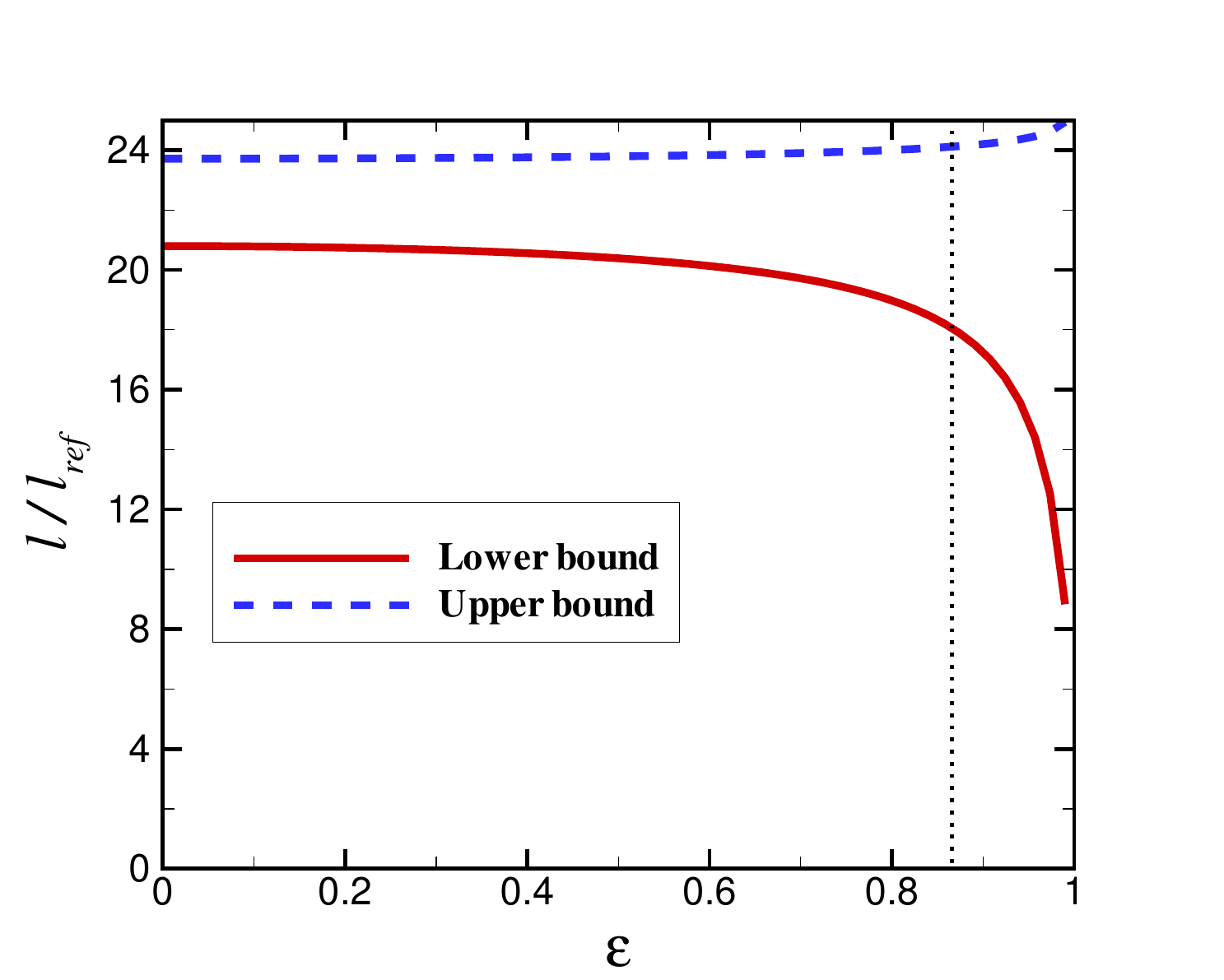} \label{fig:ltotal_ellipse_inclinee_epsilon}
 }
    \caption{Scaled effective lengths and end correction for an elliptical perforation with a $\theta=60$\textdegree~tilting angle. $\epss$ is the eccentricity of the elliptical opening. The dotted vertical line 
    corresponds to the eccentricity for which the cross-section of the perforation is circular. }
  \label{fig:l_incline_epsilon}
 \end{figure}

The eccentricity of the opening strongly affects the Rayleigh conductivity when compared to the circular case corresponding to $\epss=0$. As the eccentricity tends towards 1, in other words the perforation is changed into a slit of finite width $2b$ but infinite length $2a$, the Rayleigh conductivity tends to infinity and a divergence between lower and upper bounds can be observed. This divergence is all the more visible on the effective length and end correction plots (see figure \ref{fig:l_incline_epsilon}). 
 
However, it is worth noting that in the opposite case where the perforation is changed into a slit of zero width and finite length (\ie $a$ is kept constant and  $b$ tends to zero), then the Rayleigh conductivity tends to zero and the lower and upper bounds converge.
  
All these results show that geometrical features greatly affect the Rayleigh conductivities or end corrections of perforations, and thus cannot be discarded if  the effects of this parameter on the acoustic properties of the plate are to be taken into account accurately. Of course, this is only true if the potential component of the plate impedance (\ie the impedance due to potential flow within the perforations, as modelled in the present study)  is significant, which is not the case for liners in grazing flow ducts, for instance. However, the potential form of the Rayleigh conductivity is the basis   for models used in bias flow configurations, and an accurate description of it can be of great help, as will be shown in the next section.

\subsubsection{Eldredge \etal\ numerical simulations}
The influence of the geometrical features can also be illustrated from the  numerical simulations of \citet{2007_AIAA_Eldredge} who performed incompressible  large-eddy simulations of the acoustic response of a turbulent flow through a multiperforated liner. The perforations had a circular bore of diameter $2r=5$mm,   a tilting angle $\theta = 60^o$, and a thickness of $h = 10$mm.  Thus, at the bottom and top of the plate, the openings of the perforations   were elliptical with a major axis of $10$mm and a minor axis of $5$mm. However, the cross section along the tilted axis was circular with a cross-section of $\pi r^2$. This case corresponds exactly to the configuration considered in expression (\ref{eq:ellipseinc_bis}).

We recall that for this configuration, if the reference section is chosen to be $\pi r^2$, the effective height to be considered is $\widetilde{h}=h / \cos \theta$.  This is the value used by Eldredge \etal\ for calculating the Rayleigh conductivity of the hole, but they applied it with a modified version of Howe's model in which the no-flow Rayleigh conductivity was in fact the upper bound of the untilted cylinder bounds (\ref{eq:bounds-untiltedcylinder}):
 \begin{equation}\label{eq:Eldredge}
 K_R^{\textrm{Eld}}=\frac{ \pi r^2}{ \widetilde{h} + \pi r / 2}
 \end{equation} 

 Their model's results are similar to those of their simulation, but they obtained a better agreement by taking the effective height to be  $0.75 h/\cos \theta$. They explained that this decrease of the effective height may be due to the fact that the radius of the aperture is effectively  larger than $r$ because of the elliptical intersection of the aperture with the top and bottom of the plate.

Table \ref{tab:geom-Eldredge} summarizes the values of the no-flow Rayleigh conductivity proposed by Eldredge \etal\ and the analytical bounds obtained using the expressions established in section \ref{section:cond}. The Eldredge \etal\ value obtained with $\widetilde{h}=0.75 h/\cos \theta$ is within  the bounds for the tilted perforation with an elliptical opening given by (\ref{eq:ellipseinc_bis}), and is quite close to the mean value of the bounds  $K_R = (K_R^- + K_R^+) /2 = 11,375.10^{-4}$. On the contrary, their first value obtained with $\widetilde{h}= h/\cos \theta$, is lower than the lower bound for   a tilted perforation with an elliptical opening and greater than the upper bound for a tilted perforation with a circular opening of radius $r$ given by (\ref{eq:cylinc}). This points out the importance of the tilted angle in the computation of the Rayleigh conductivity. 

The coefficient $0.75 / \cos \theta$ introduced by Eldredge \etal\ corrects empirically the difference between the bounds for an untilted cylinder (\ref{eq:bounds-untiltedcylinder}) and for a tilted perforation with elliptical opening (\ref{eq:ellipseinc_bis}).  The bounds (in \ref{eq:ellipseinc_bis}) therefore improve the accuracy of the Rayleigh conductivity approximation, by taking into account the dependency on the parameters $a$ and $b$ of the ellipse, as well as on the tilted angle $\theta$.

\begin{table}
  \begin{center}
\def~{\hphantom{0}}
  \begin{tabular}{ccccc} 
 &  Tilted perforation, & Tilted perforation, &Expression (\ref{eq:Eldredge}) & Expression (\ref{eq:Eldredge})  \\
& circular opening  & elliptical opening& with $\widetilde{h}=h/\cos \theta$ & with $\widetilde{h}=0.75h/\cos \theta$\\[3pt]
   $K_R^-$&    $4.44\,10^{-4}$ & $8.57\,10^{-4}$ &\multirow{2}{*}{$8.21\,10^{-4}$}  & \multirow{2}{*}{$10.00\,10^{-4}$}\\[2pt]
   $K_R^+$&       $5.73\,10^{-4}$ & $14.18\,10^{-4}$ &  &\\
  \end{tabular}
  \caption{Bounds of the Rayleigh conductivity obtained from expressions established in section \ref{section:cond}, compared with the values used in the paper by Eldredge \etal\ .}
  \label{tab:geom-Eldredge}
  \end{center} 
\end{table}

These results confirm that there is an advantage in taking the geometry of the aperture into account correctly. Of course, the question of the modification of the no-flow Rayleigh conductivity in the presence of a bias flow is still open. Strategies used until now, which consist in modifying this value by a flow-dependent function  (either as done by \citet{2000_AIAAJ_Jing} or in a slightly different way by \citet{2005_JFS_Howe}) may still be relevant.


\section{Effective acoustic compliance of a low-porosity perforated plate}\label{section:reseau}

In this section, we focus on a method for determining the asymptotic expansion of the reflection and transmission coefficients of the aforementioned low-porosity perforated plate up to order 2 relative to the small parameter characterizing the ratio of a characteristic size $d$ of the perforation to the spacing $L$ between two neighbouring perforations. This expansion has been recently given in \citet{bendali2012}. However, contrary to the approach in this reference, which was based on complex two-scale matched asymptotic expansions and the consideration of a grating of multipoles, the present one deals with these coefficients in a direct way without resorting to the expansion of the whole wave. Moreover it makes it possible to elucidate the role played by some relevant physical parameters as a weighted flux in which these coefficients are expressed. These new features are certainly more compatible for use in a fluid  mechanics context. 

The method is an adaptation of an approach previously used by \citet{1973_JFM_Leppington} for dealing with infinitely thin perforated plates. Actually, the method devised by Leppington and Levine is extended here in several directions. The plate is no longer assumed to be infinitely thin and the perforations can be of a rather arbitrary shape. As mentioned above, it is worth recalling that the perforated plates under consideration in this study are assumed to be of low-porosity, \ie such that $\delta = d/L \ll 1$, and  to behave relative to the incident acoustic wave as a homogeneous surface, \ie such that the spacing $L$ is less that a half-wavelength, $L< \lambda /2$. 

The main mathematical tools are asymptotic expansions, integral equations, and the lattice sum theory for the Helmholtz equation \citep{linton2010}. We end this section by briefly recalling results from \citet{bendali2012} on the derivation of effective compliances for the plate from the asymptotic expansion of the reflection and transmission coefficients.

\subsection{Scattering problem and Floquet's modes}\label{subsection:2.1}

Consider the configuration shown in figure \ref{fig:plaque-schema} of an acoustic plane wave incident upon a perforated plate. As is well-known in lattice theory \citep{nedelec1991, linton2010}, the total wave satisfies the following quasi-periodic conditions resulting from the periodicity properties of the array and the fact that translating the incident wave by a vector $\bxi_{m}$, linked to the periodicity properties of the lattice (see (\ref{eq:latticeproperties})), simply shifts the phase according to the following formula 
\begin{equation}
p(\bx+\bxi_{m})=e^{\ci\bbeta\bcdot\bxi_{m}}p(\bx),\label{eq:quasiperiodic}
\end{equation}
where $\bbeta\cdot\bxi_{m}$ is the scalar product of the Bloch vector $\bbeta=\kappa\btau'$ and $\bxi_{m}$. In the following, when saying that $p$ satisfies periodic boundary conditions, we mean that $p$ fulfils condition (\ref{eq:quasiperiodic}). This of course induces the same property on any derivative of $p$. It is therefore sufficient to restrict the determination of the wave $p$ to the unit cell displayed in Fig. \ref{fig:domain}.

This determination requires a radiation condition on the reflected wave $\big(p-p_{i}\big)\vert_{\sU_{+}}$ and the transmitted one $p\vert_{\sU-}$, which can be set in terms of a series of Floquet modes \citep[for instance, see][]{bendali2012, linton2010, nedelec1991}.

As a result, the wave reflected and transmitted by the multiperforated plate can be specified in the unit cell by solving the following boundary-value problem: 
\begin{equation}
\left\{ \begin{array}{ll}
\nabla^{2}p+\kappa^{2}p=0 & \mbox{ in }\sU\\
\partial_{\bn}p=0 & \mbox{ on the plate},
\end{array}\right.\label{pb:scattering}
\end{equation}
together with outgoing radiation conditions and periodicity properties (\ref{eq:quasiperiodic}) on the fictitious boundaries $\partial\sU_{\pm}$, for $p$ and its derivatives up to order 1.


\subsection{A domain decomposition approach}

\label{subsection:2.2}

A boundary integral equation formulation is used on $\sU_{\pm}$. It should be noted that the plate is assumed to be located between $x_{3}=-h/2$ and $x_{3}=h/2$ and that a local coordinate system $x_3^\pm = x_3 \pm h/2$ is used on $\sU_{\pm}$ so that $x_3^+ = 0$ and $x_3^-=0$ correspond to the upper and lower side of the plate respectively.

In \citet{1973_JFM_Leppington}, the main tool for dealing with the wave outside the perforations is the quasi-periodic Green's kernel for the Helmholtz equation $G_{\Lambda}$: 
\begin{equation}
G_{\Lambda}(\bx,\by)=\sum_{m}e^{\ci\bbeta\bcdot\bxi_{m}}G(\left(\bx'-\by',\vert x_{3}\vert\right),\bxi_{m}),\label{eq:quasiperiodickernel}
\end{equation}
with $\bx=(\bx',x_{3})$, $\by=(\by',0)$. The underlying kernel 
\begin{equation}
G(\bx,\bxi)=\frac{e^{\ci\kappa\absxxi}}{4\pi\absxxi},
\end{equation}
is the standard Green's function yielding the outgoing solution to the Helmholtz
equation. The solution of the scattering problem when the plate is
considered as perfectly reflecting can be explicitly written
as: 
\begin{equation}
\left\{ \begin{array}{ll}
p_{0}(\bx_{+})=p_{i}\left(\bx'_{+},x_{3}^{+}\right)+p_{i}\left(\bx'_{+},-x_{3}^{+}\right)\\
p_{0}(\bx_{-})=0 & .
\end{array}\right.\label{eq:u0}
\end{equation}

Then, the solution to  boundary-value problem (\ref{pb:Helmholtz}) can be expressed, outside the perforation, in terms of a single
integral equation written in local coordinates: 
\begin{equation}
p_{\pm}(\bx_{\pm})=p_{0}(\bx_{\pm})\mp\int_{\sD_{\pm}}2\: G_{\Lambda}(\bx_{\pm},\by')\partialtroispm(\by',0)\: \dif\by'\label{eq:integralequation}
\end{equation}

Indeed, the quasi-periodic Green kernel $G_{\Lambda}$ takes into
account all the following features: the Helmholtz equation, the quasi-periodic
conditions and the outgoing radiation conditions. Moreover, the actual
kernel here is $2G_{\Lambda}$ which, from the method of images, ensures
that $\partial_{x_{3}}p_{\pm}$ obtained as limiting value from the
integral representation formula is 0 outside the openings of the perforation
and equal to their value assigned under the integral inside them.
As a result, the acoustic pressure within the hole can hence be obtained
by solving the following system: 
\begin{equation}
\left\{ \begin{array}{ll}
\begin{array}{ll}
\nabla^{2}p+\kappa^{2}p=0 & \mbox{ on }\Omega\\
\partial_{\bn}p=0 & \mbox{ on }\Sigma
\end{array}\\
\begin{array}{ll}
p_{\pm}(\bx_{\pm})=p_{0}(\bx_{\pm})\mp\ds\int_{\sD_{\pm}}2\: G_{\Lambda}(\bx_{\pm},\by)\partialtroispm(\by',0)\: \dif\by'.\end{array}
\end{array}\right.\label{pb:DDM}
\end{equation}
combined with transmission conditions between the hole and the integral
formulation within the upper and lower parts of the unit cell: 
\begin{equation}
\left\{ \begin{array}{ll}
p=p_{+}\mbox{ on }\sD_{+},\\
p=p_{-}\mbox{ on }\sD_{-},
\end{array}\right.\qquad\mbox{ and }\qquad\left\{ \begin{array}{ll}
\partialtrois=\partialtrois_+\mbox{ on }\sD_{+},\\
\partialtrois=\partialtrois_{-}\mbox{ on }\sD_{-}.
\end{array}\right.\label{eq:transmissionconditions}
\end{equation}

It is worth noting that the above transmission conditions are expressed
in terms of functions defined on the same sets instead of functions
depending on the same variables. It is in this context that the domain
decomposition approach can handle efficiently the incident and reflected
waves on one hand and the transmitted wave on the other hand.

\subsection{Second-order asymptotic expansions}\label{subsection:LLM}

Since $L<\lambda/2$, dual, or spectral, representations of the Green's kernel $G_{\Lambda}$ for $|x_{3}|\gg1$ \citep[see, \eg][formula (2.9)]{linton2010} make it possible to obtain the reflection and the transmission coefficients, respectively denoted $R_{\delta}$ and $T_{\delta}$ to underline their dependence on $\delta = d/L$, through the following decomposition of the wave into a propagative and an evanescent part: 
\begin{equation}
\begin{array}{ll}
p(\bx_{+}',x_{3}^{+})=e^{\ci\kappa\btau'\bcdot\bx'_{+}}\left(e^{-\ci\kappa\cos{\Phi}x_{3}^{+}}+R_{\delta}e^{\ci\kappa\cos{\Phi}x_{3}^{+}}\right)+\mbox{evan. modes} & \mbox{ for }x_{3}^{+}>0,\\
p(\bx_{-}',x_{3}^{-})=e^{\ci\kappa\btau'\bcdot\bx'_{-}}\: T_{\delta}e^{-\ci\kappa\cos{\Phi}x_{3}^{-}}+\mbox{evan. modes} & \mbox{ for }x_{3}^{-}<0.
\end{array}\hspace{-0.5em}
\end{equation}

These decompositions are obtained as follows. A Poisson summation
formula is used to express the quasi-periodic Green's kernel in terms
of propagative and evanescent modes \citep{linton2010}: 
\begin{equation}
G_{\Lambda}(\bx,\by)=\frac{\ci}{2A}\sum_{m}\frac{1}{\gamma_{m}}e^{\pm\ci\gamma_{m}x_{3}}\, e^{\ci\bbeta_{m}\bcdot(\bx'-\by')},
\end{equation}
where $\gamma_{m}\ds=\sqrt{\kappa^{2}-|\bbeta_{m}|^{2}}$ and $\bbeta_{m}=\bbeta+2\pi\left(m_{1}\bxi_{1}^{*}+m_{2}\bxi_{2}^{*}\right)$,
with $\left(\bxi_{1}^{*},\bxi_{2}^{*}\right)$ the dual basis of
$\left(\bxi_{1},\bxi_{2}\right)$. Under condition (\ref{eq:2.6})
on the lattice, only the fundamental mode is propagating, so that
\begin{equation}
G_{\Lambda}(\bx,\by)=\frac{\ci}{2\kappa A\cos{\Phi}}\: e^{\ci\kappa\btau'\bcdot(\bx'-\by')}+\mbox{ evanescent modes}.
\end{equation}

The reflection and transmission coefficients are hence 
\begin{equation}
R_{\delta}=1-\frac{\ci Q_{\delta}^{+}}{\kappa A\cos{\Phi}}\:\mbox{ and }\: T_{\delta}=\frac{\ci Q_{\delta}^{-}}{\kappa A\cos{\Phi}},\label{eq:4.X3}
\end{equation}
where $Q_{\delta}^{\pm}$ are the following weighted fluxes~: 
\begin{equation}
Q_{\delta}^{\pm}=\ds\int_{\sD_{\pm}}e^{-\ci\kappa\btau'\bcdot\by'}\partial_{x_{3}}p_{\pm}(\by',0)\:\dif\by'.\label{eq:4.X4}
\end{equation}

The objective in this section is to adapt the approach used by \citet{1973_JFM_Leppington}
to get a second-order asymptotic expansion for $R_{\delta}$ and $T_{\delta}$
in powers of $\delta$ for general perforations and plates.

Equations (\ref{pb:DDM}) and (\ref{eq:transmissionconditions}) make
it clear that the determination of $p$, $p_{+}$ and $p_{-}$ only
requires solving a problem set solely on the hole $\Omega$. To suppress
the dependence of the geometry on $\delta$, it is convenient to introduce
the scaled variable 
\begin{equation}
\bX=\bx/\delta,\label{eq:4.X5}
\end{equation}
also called a fast (or inner) variable in asymptotic expansions theory.
In this way, $\Omega$, $\Sigma$ and $\mathsf{D_{\pm}}$ are parameterized
by fixed sets respectively denoted $\homega$, $\hsigma$ and $\hsDpm$,
meaning that they are independent of $\delta$. Domain $\homega$
can hence be viewed as an isolated perforation in a plate extending
to infinity, see \citet{tuck1975}, \citet{Howe_book} (where
this feature is introduced in a quite implicit way) and \citet{bendali2012}
where it is dealt with by means of a matched asymptotic expansion
technique.

The expansion of $R_{\delta}$ and $T_{\delta}$ will be obtained
through the second-order asymptotic expansion of the function $\Pi_{\delta}$
defined on $\homega$ through the scaling (\ref{eq:4.X5}) 
\begin{equation}
p(\delta\bX)=\Pi_{\delta}(\bX)=\Pi^{(0)}(\bX)+\delta\Pi^{(1)}(\bX)+\circ\left(\delta\right).\label{eq:4.X6}
\end{equation}

Clearly, introducing expansion (\ref{eq:4.X6}) into the Helmholtz
equation and the boundary condition on $\Sigma$ yields 
\begin{equation}
\left\{ \begin{array}{ll}
\nabla_{\bX}^{2}\Pi^{(j)}=0 & \mbox{ in }\homega\\
\partial_{\bn}\Pi^{(j)}=0 & \mbox{ on }\hsigma,
\end{array}\right.\label{syst:1bis}
\end{equation}
for $j=0,1$, with $\bn$ the unit normal to $\hsigma$ directed outwards
$\homega$. The expansion of the integral equations set on respectively
$\sD_{+}$ and $\sD_{-}$ requires using the theory of lattice sums
for the Helmholtz equation, which gives \citep[page 653]{linton2010}
\begin{equation}
G_{\Lambda}(\delta\bX,\delta\bY)=\frac{1}{4\pi\delta\vert\bX'-\bY'\vert}+s_{0}+\circ\left(\delta\right),
\end{equation}
where $s_{0}$ has an explicit expression in terms of a Schl\"omilch
series, see \citet{linton2009}. As a result, expanding 
\begin{equation}
p_{0}(\delta\bX'_{+},0)=2e^{\ci\kappa\btau'\bcdot\bX'_{+}}=2+2\delta\ci\kappa\btau'\bcdot\bX'_{+}+\circ(\delta),\label{eq:3}
\end{equation}
and expressing the integral equations of (\ref{pb:DDM}) in terms
of the fast variables, gives 
\begin{systeme}\label{syst:pizero}
\ds \pizero_+(\bX'_+) + \int_{\hsD_+} \frac{\partial_{X_3} \pizero_+
(\bY')}{2 \pi \vert \bX'_+-\bY' \vert} \: \dif \bY' =2,\\
 \ds \pizero_-(\bX'_-) - \int_{\hsD_-} \frac{\partial_{X_3}
\pizero_- (\bY')}{2 \pi \vert \bX'_--\bY' \vert} \: \dif \bY' =0,
\end{systeme}
 \begin{systeme}\label{syst:piun} 
 \piun_+(\bX'_+) + \int_{\hsD_+} \frac{\partial_{X_3} \piun_+(\bY') }{2 \pi
\vert \bX'_+-\bY' \vert} \: \dif \bY' = -2 s_0 \int_{\hsD_\pm}
\partial_{X_3} \pizero \dif \bY' + 2 \ci \kappa \btau'\bcdot \bX'_+, \\
\ds \piun_-(\bX'_-) - \int_{\hsD_-} \frac{\partial_{X_3} \piun_-(\bY')}{2 \pi \vert \bX'_--\bY' \vert} \: \dif \bY' =\phantom{-}
2 s_0 \int_{\hsD_\pm} \partial_{X_3} \pizero \dif \bY',
\end{systeme}
 where $\Pi_{\pm}^{(j)}$ stand for the respective expressions
of $\Pi^{(j)}$ on $\hsD_{\pm}$ in terms of the variables $\bX'_{\pm}=\bx'_{\pm}/\delta$.
Define $\pizero$ in $\homega$ so that it can be extended by 
\begin{equation}
\pizero_{+}(\bX'_{+})=2-\int_{\hsD_{+}}\frac{\partial_{X_{3}}\pizero_{+}(\bY')}{2\pi\vert\bX'_{+}-\bY'\vert}\:\dif\bY'\mbox{ and }\pizero_{-}(\bX'_{-})=\int_{\hsD_{-}}\frac{\partial_{X_{3}}\pizero_{-}(\bY')}{2\pi\vert\bX'_{-}-\bY'\vert}\:\dif\bY',
\end{equation}
for $\pm X_{3}^\pm >0$, where $ X_{3}^\pm $ are the component normal to the plate relatively to the coordinates $\bX_\pm$ respectively. Actually, equations (\ref{syst:1bis})
and (\ref{syst:pizero}) mean that $\pizero$ is the solution of the
following boundary-value problem 
\begin{equation}\label{eq:4.X14}
\left\{ \begin{array}{ll}
&\nabla_{\bX}^{2}\Pi^{(0)}=0 \qquad \mbox{ in }\homega,\\
&\partial_{\bn}\Pi^{(0)}=0  \:  \:  \qquad \mbox{ on }\hsigma,\\
& \lim_{X_{3}\rightarrow+\infty}\pizero_{+}=2, \\
& \lim_{X_{3}\rightarrow-\infty}\pizero_{-}=0.
\end{array} \right.
\end{equation}

This shows that 
\begin{equation}
\hKR=\frac{1}{2}\int_{\hsD_{\pm}}\partial_{X_{3}}\pizero\dif\bY'.\label{eq:4.X15}
\end{equation}
is nothing else than the Rayleigh conductivity of the isolated perforation
in the scaled variables (\ref{eq:4.X5}). Coming back to the initial
variables, we can express $\hKR$ by means of the usual Rayleigh conductivity
defined in (\ref{eq:Ray-cond-p-adim}) 
\begin{equation}
\hKR=K_{R}/\delta.\label{eq:4.X16}
\end{equation}

In a similar way, $\piun$ can be expressed from the solution of the
following problem 
\begin{equation}
\left\{ \begin{array}{ll}\label{eq:4.X17}
&\nabla_{\bX}^{2}\Pi_{l}^{\bX}=0 \qquad \mbox{ in }\homega\\
&\partial_{\bn}\Pi_{l}^{\bX}=0  \:  \:  \qquad \mbox{ on }\hsigma,\\
&\lim_{X_{3}\rightarrow+\infty}\left(\Pi_{l}^{\bX}(\bX)-\bX'_{l}\right)=2,\\
&\lim_{X_{3}\rightarrow-\infty}\Pi_{l}^{\bX}(\bX)=0,
\end{array}\right.
\end{equation}
for $l=1,2$ and the previous function $\pizero$ through 
\begin{equation}
\piun(\bX)=-4s_{0}\:\hKR\:\pizero(\bX)+2\ci\kappa\btau'\bcdot\bPi^{\bX}(\bX),\label{eq:4.X18}
\end{equation}
where $\bPi^{\bX}$ is the vector function whose components are respectively
$\Pi_{1}^{\bX}$, $\Pi_{2}^{\bX}$ and $0$.

The fluxes $Q_{\delta}^{\pm}$ defined in equation (\ref{eq:4.X4})
now become, up to some $\circ(\delta^{2})$ terms, 
\begin{equation}\label{eq:4.X19}
\begin{array}{ll}
Q_{\delta}^{\pm} & =\ds\delta\int_{\hsD_{\pm}}\exp\left(-\ci\kappa\delta\btau'\bcdot\bY'\right)\left(\partial_{X_{3}}\pizero_{\pm}+\delta\partial_{X_{3}}\piun_{\pm}\right)(\bY'_{\pm},0)\:\dif\bY'+\circ(\delta^{2})\\
 & \ds=2\delta\hKR+\delta^{2}\int_{\hsD_{\pm}}\left(\partial_{X_{3}}\piun_{\pm}-\ci\kappa\btau'\bcdot\bY'\partial_{X_{3}}\pizero_{\pm}\right)\:\dif\bY',
\end{array}
\end{equation}
which, thanks to (\ref{eq:4.X18}), gives 
\begin{equation}\label{eq:4.X20}
Q_{\delta}^{\pm}=2\delta\hKR+\delta^{2}\left(2\ci\kappa\btau'\cdot\int_{\hsD_{\pm}}\left(\partial_{X_{3}}\bPi_{\pm}^{\bX}-\frac{1}{2}\bY'\partial_{X_{3}}\pizero_{\pm}\right)\:\dif\bY'-8s_{0}\hKR^{2}\right).
\end{equation}
A reciprocity property \citep[see][page 16]{bendali2012} yields
 \begin{systeme}\label{eq:4.X21} 
& \ds \int_{\hsD_+} \left( \partial_{X_3} \bPi^\bX_+ - \frac{1}{2}\bY'+ \partial_{X_3} \pizero_+ \right) \dif \bY'_+ = 0 \\
& \ds \int_{\hsD_-} \left( \partial_{X_3} \bPi^\bX_- - \frac{1}{2}\bY'- \partial_{X_3} \pizero_- \right) \dif \bY'_- = -\widehat{\bmu_n}, 
\end{systeme} with 
\begin{equation}
\widehat{\bmu_{n}}=(\bc_{+}-\bc_{-})\hKR-\frac{1}{2}\int_{\hsigma}\pizero\bn'\:\dif s,
\end{equation}
where $\bc_{\pm}=\bX-\bX_{\pm}$ take into account the difference
of the centres of phases chosen for the reflected and transmitted
waves and $\bn'$ is the projection of the unit normal $\bn$ to $\hsigma$
on the plane of the plate. For an axisymmetric perforation, $\widehat{\bmu_{n}}$
reduces to 
\begin{equation}
\widehat{\bmu_{n}}=\frac{h}{\delta}\:\hKR\:\be_{3}=\frac{h\: K_{R}}{\delta^{2}}\:\be_{3},\label{eq:4.X25}
\end{equation}
where $\be_{3}$ is the unit vector along the $x_{3}$-axis.

As a result, the reflection and the transmission coefficients have
the following asymptotic expansions, as obtained in \citep{bendali2012}
through a more involved procedure
expansions 
\begin{systeme} \label{eq:4.X31} 
\ds & R_\delta = \ds 1 +\frac{ \delta \: \hKR }{ A} \frac{ 2 }{\ci \kappa
\cos{\Phi}} - \frac{ \delta^2 \: 4 \: \hKR^2 \: s_0 }{A }\frac{ 2 }{\ci \kappa \cos{\Phi}} + \circ(\delta^2) \\[1em] 
\ds & T_\delta =\phantom{1} -\ds \frac{ \delta \: \hKR }{ A} \frac{ 2 }{\ci \kappa \cos{\Phi}} + \frac{
\ds \delta^2 \: \left(4 \: \hKR^2 \: s_0 + \ci \kappa \btau' \cdot \widehat{\bmu_n} \right)}{ A } \ds \frac{ 2 }{\ci \kappa \cos{\Phi}} + \circ(\delta^2). 
\end{systeme} 

It is worth mentioning that the term $\btau'\bcdot\widehat{\bmu_{n}}$ is
equal to zero for an axisymmetric perforation since then $\btau'$
has no component along the $x_{3}$-axis.

Finally, in terms of the Rayleigh conductivity in the initial variables,
we get \begin{systeme} \label{eq:4.X32} 
\ds & R_\delta = \ds 1
-\frac{2 \ci }{\kappa A \cos{\Phi}} K_R + \frac{8 \ci
s_0 }{ \kappa A \cos{\Phi}} K_R^2 + \circ(\delta^2) \\[1em]
\ds & T_\delta = \phantom{1-} \ds \frac{2 \ci }{\kappa A
\cos{\Phi}} K_R - \frac{ 2 \ci s_0 }{\kappa A \cos{\Phi}}\left(
4 K_R^2 + \ci \kappa \btau' \cdot \bmu_n \right) + \circ(\delta^2).
\end{systeme}

\subsection{Effective compliance of the perforated plate}\label{sec:4.4}

The effective compliance $K$ of the plate can hence be defined through
the following relationships 
\begin{equation}
\partial_{x_{3}}p_{+}=\partial_{x_{3}}p_{-}=K(p_{+}-p_{-}).\label{eq:4.X28}
\end{equation}

The reflection and transmission coefficients and the compliance are
then linked by 
\begin{equation}
R+T=1,\quad K=\frac{\ci\kappa\cos{\Phi}}{2}\left(1-\frac{1}{R}\right),\quad R=1/\left(1-\frac{2K}{\ci\kappa\cos{\Phi}}\right)\label{eq:4.X29}
\end{equation}

Moreover, if, as here, the governing equations are those of linear
acoustics, which means that no damping mechanism is involved and that
there is no acoustic absorption by the perforated plate, the imaginary
part $\Im(K)$ of the effective compliance is zero and the conservation
of energy is expressed through the following relation 
\begin{equation}
1-|R|^{2}-|T|^{2}=0.\label{eq:4.X30}
\end{equation}

Following a general procedure given in  \citet{bendali2012}, it is then possible
from expansion (\ref{eq:4.X32}) to design a compliance which is both
consistent with the approximation of the effective reflection and
transmission coefficients at leading order $\bigcirc(\delta)$, and
ensures the conservation of acoustic energy: 
\begin{equation}\label{eq:4.X35}
K_{\delta}^{(1)}=\frac{K_{R}}{A}.
\end{equation}

This is the classical expression of the effective compliance of the
perforated plate \citep{2003_JFM_Eldredge, 1990_JFM_Hughes}. The
reflection coefficient corresponding to this compliance condition
\begin{equation}\label{eq:4.X36}
R_{\delta}^{(1)}=\frac{1}{\ds1-\frac{2}{\ci\kappa\cos{\Phi}}\frac{K_{R}}{A}}
\end{equation}
does not lead to an amplified reflected wave while approximating the
actual reflection coefficient up to a term in $\circ(\delta)$. 

In general, the same procedure cannot be repeated to get a second-order
compliance condition since the term $\btau'\cdot\widehat{\bmu_{n}}$
in equations (\ref{eq:4.X31}) or (\ref{eq:4.X32}) prevents the condition
$R+T=1$ from being satisfied. However, this can be done when $\btau'\cdot\widehat{\bmu_{n}}=0$,
for instance when the perforation is axisymmetric or when the thickness
of the plate can be neglected. 

A very far-reaching, by no means obvious, property of the lattice
sums for the Helmholtz equation \citep[see][page 656]{linton2010}
ensures that 
\begin{equation}\label{eq:4.X41}
\Im\left(2s_{0}\right)=\frac{1}{\kappa\cos{\Phi}},
\end{equation}
and then makes it possible to design a compliance at leading order
$\bigcirc(\delta^{2})$ 
\begin{equation}\label{eq:4.X42}
K_{\delta}^{(2)}=\frac{K_{R}}{A}\left(1-4K_{R}\: \Re(s_{0})\right),
\end{equation}
\citep[cf.][]{bendali2012}. The coefficient $\Re(s_{0})$ adds a correction
depending on the shape of the array of the perforations rather than
just averaging by the area $A$ of the lattice unit cell. The corresponding
reflection coefficient 
\begin{equation}
R_{\delta}^{(2)}=\left(1+\frac{2\ci}{\kappa\cos{\Phi}}\ds\frac{K_{R}}{A}\left(1-4K_{R}\:\Re(s_{0})\right)\right)^{-1}\label{eq:4.X43}
\end{equation}
does not exhibit any unrealistic behaviour. \citet{1973_JFM_Leppington}
gave a similar expression for a rectangular lattice of elliptical
holes (see their equation (2.27)). Our approximation
(\ref{eq:4.X43}) of the reflection coefficient is an extension of
this formula to an arbitrary shape of the lattice and to any geometry
of the perforations, if either the plate is infinitely thin or the
perforations are axisymmetric. Otherwise, second-order terms in $\bigcirc(\delta^{2})$ are missed.


\subsection{Numerical results for a typical plate}

Some numerical experiments can now be conducted to illustrate the influence of the geometry of the perforations and the shape of the lattice for a realistic configuration of the perforated plate. A $2$mm perforated plate of $1.96\%$ porosity is considered, with either a rectangular or staggered lattice of holes. As in section \ref{sec:generic-perforated-plate}, the characteristic size of the holes is $r=0.225$mm. The incidence of the incoming acoustic wave is $\Phi=45^{o}$. The dimensions of the lattice are $\norm{\bxi_{1}}=3$mm and $\norm{\bxi_{2}}=2.7$mm in the rectangular configuration, $\norm{\bxi_{1}}=3$mm and $\norm{\bxi_{2}}=\sqrt{(1.5)^{2}+(2.7)^{2}}\approx3.1$mm in the staggered one. Thus, the surface of the lattice cell is the same in both configurations, i.e., $A=8.1$mm$^{2}$. The assumption that $\delta=r/L$ is a small parameter is satisfied.

Here, we focus on two kinds of geometries for the holes: the untilted cylinder of radius $r$ and the tilted perforation of elliptical opening with circular cross-section of radius $r$ and a $60^{o}$ tilt angle. The mean value of the bounds (\ref{eq:bounds-untiltedcylinder}) and (\ref{eq:ellipseinc_bis}) are taken as the respective Rayleigh conductivities of these perforations. As discussed in section \ref{sec:4.4}, the reflection coefficient can be computed at first and second order in $\delta$ by equations (\ref{eq:4.X36}) and (\ref{eq:4.X43}), respectively.

\begin{figure}
\centering 
\subfloat[Amplitude]{\includegraphics[width=0.45\linewidth]{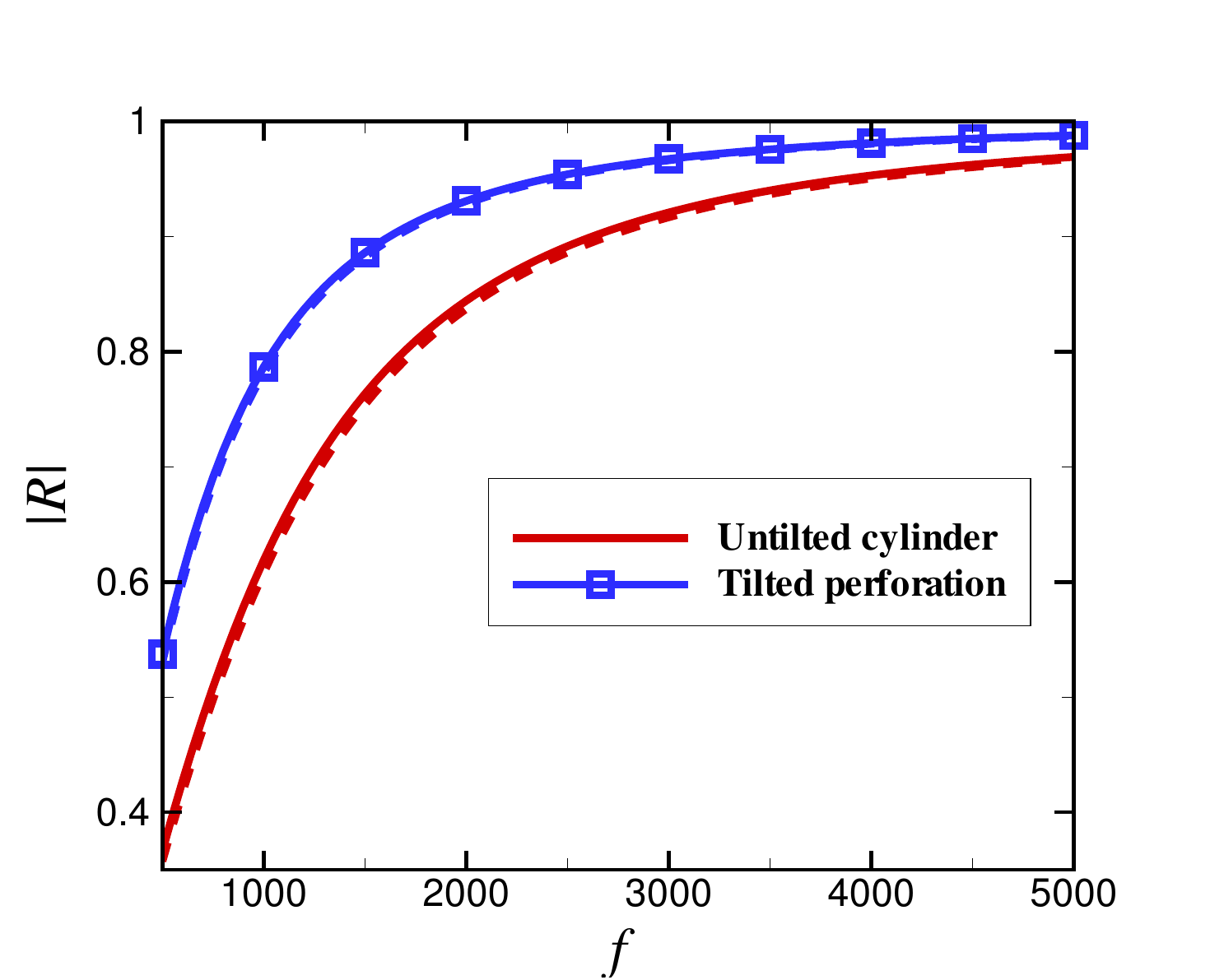} \label{fig:absR} 
}
\subfloat[Phase (in degrees)]{\includegraphics[width=0.45\linewidth]{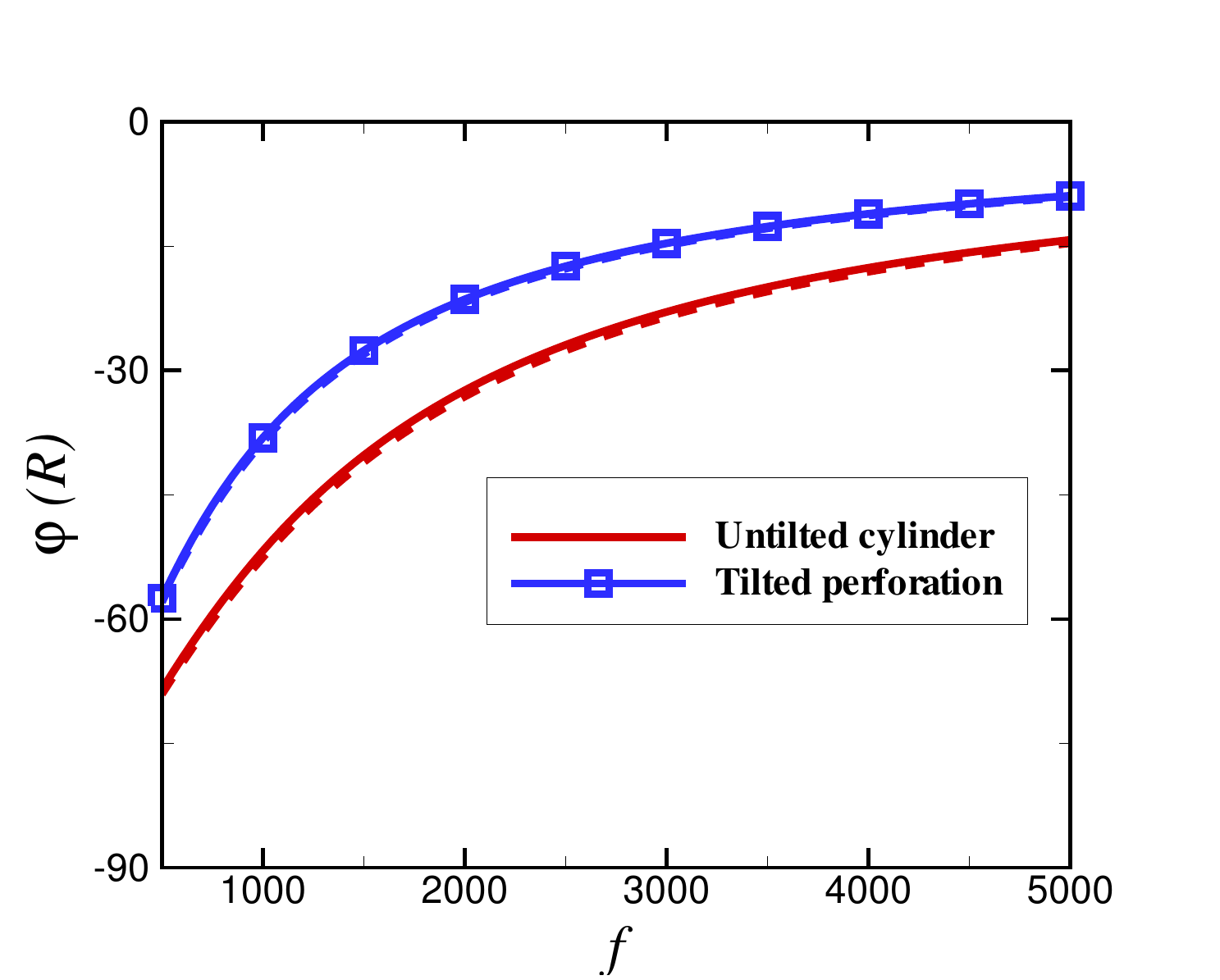} \label{fig:phiR} 
}
\caption{Reflection coefficient of the perforated plate with a rectangular
lattice of untilted cylindrical perforations (red lines without symbols)
or tilted perforations with an elliptical opening and circular cross-section (blue lines with
symbols). First order approximations are plotted in solid lines, second
order ones in dashed lines.}
\label{fig:reseaudroit} 
\end{figure}

Results in the frequency range $[500,5000]$Hz are plotted in figure
\ref{fig:reseaudroit}, for the rectangular lattice. The geometry
of the perforations has a strong effect on the reflection coefficient,
which confirms the importance of using the expressions for the Rayleigh
conductivity obtained in section \ref{section:cond}. On the contrary,
first and second order approximations are almost indistinguishable,
which shows that the usual way of calculating the effective compliance
of the plate (i.e., by averaging the Rayleigh conductivity of the
isolated hole by the area of the lattice unit cell) is quite accurate
in this case.

\begin{figure}
\centering 
\subfloat[Amplitude]{\includegraphics[width=0.45\linewidth]{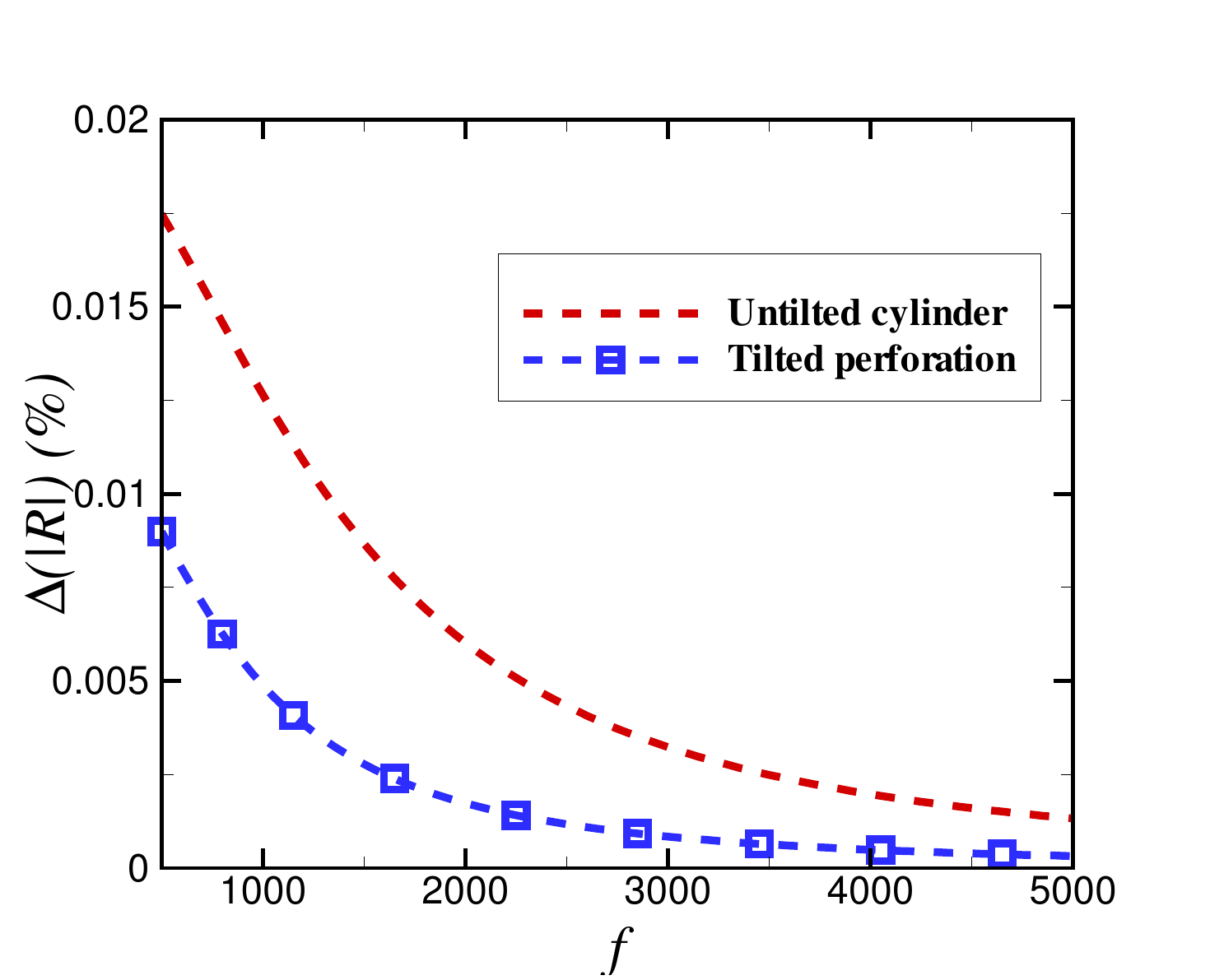} \label{fig:Delta_absR} 
}
\subfloat[Phase ]{\includegraphics[width=0.45\linewidth]{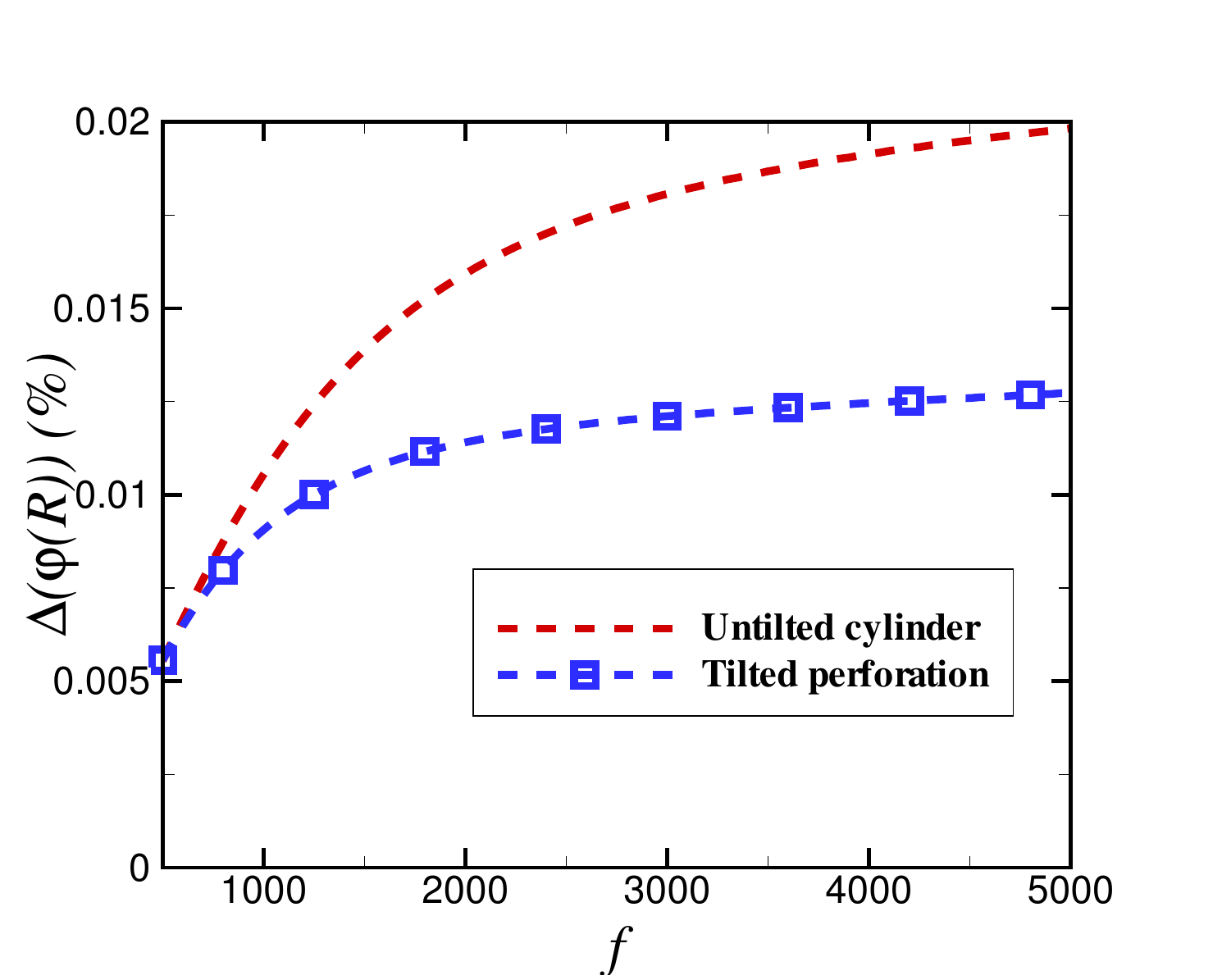} \label{fig:Delta_phiR} 
}
\caption{Difference between the second-order reflection coefficient in the
staggered and rectangular arrangement of the holes, normalized by
the rectangular arrangement value. }
\label{fig:Delta} 
\end{figure}

Consequently, it seems obvious that changing the shape of the lattice
will not affect the reflection coefficient significantly. This is
shown in figure \ref{fig:Delta}, which depicts the deviation of the
staggered array from the rectangular one expr
essed as a percentage
of the latter.

For both geometries of the perforations and whatever the frequency,
the reflection coefficients differ from each other by less than $0.02\%$.
This means that the shape of the lattice can be neglected without
reducing the quality of the models, and that considering first-order
usual approximations for the reflection and the transmission coefficients
is sufficient.

\section{Conclusion}

The acoustic properties of a low-porosity perforated plate in a compressible ideal inviscid fluid have been investigated in the absence of mean flow. This generalizes previous results done by \citet{1973_JFM_Leppington} to account for thick plates with perforations of arbitrary periodic patterns and various geometries. The study was conducted in two stages. As a first step, the mathematical framework developed in \citet{laurens2012} for obtaining accurate bounds of the Rayleigh conductivity for various geometric configurations, was explained in a fluid mechanics context and extended to cover tilted perforations with elliptical openings. Then asymptotic expansions of the effective reflection and transmission coefficients of a thick plate with a skew grating of perforations have been derived by extending the \citet{1973_JFM_Leppington} method. These results match those established in \citet{bendali2012} by means of a matched asymptotic expansion of the full wave, but the methodology to obtain them is much less complex and more suitable for a fluid mechanics context since it can be easily adapted to deal with other models of flow at the level of the perforations.

In section \ref{section:cond}, we focused on the Rayleigh conductivity of a single perforation, which can be related to the inertial end correction used in impedance models. Lord \citet{rayleigh} gave two bounds for the conductivity of a circular aperture of zero-thickness, and his work was extended to thick cylindrical  perforations by  \citet{Howe_book}. These expressions are nowadays widely used in impedance and conductivity models, even when perforations of other  geometries are considered. We have focused on tilted perforations of elliptical sections, which are of primary interest in industrial applications, especially in the  domain of engine combustors where walls are perforated by sub-millimetre tilted perforations for cooling purposes.    

In related theoretical or numerical studies, some researchers use the Rayleigh conductivity (or end correction) of an untilted cylindrical aperture with a modified plate thickness which corresponds to the length of the centreline of the tilted aperture. We have shown that exact lower and upper bounds for the  Rayleigh conductivity and end correction can be obtained by applying the Dirichlet and Kelvin variational principles. Our analysis has shown that these quantities were quite sensitive to the tilt angle and that the thickness effect requires a more precise calculation than just considering the intuitive effective value.
    
Moreover, we showed that all the expressions obtained by Lord Rayleigh and Howe for more conventional geometries can be easily derived from the present result,  by considering limiting values of the tilt angle or of the aperture eccentricity. It should be noted that though we have only considered the Rayleigh conductivity in the framework of the purely acoustic problem, \ie without any flow or non-linear effects, this solution can be used as the groundwork for deriving Rayleigh  conductivity in the presence of a bias flow through the aperture, as was done by \citet{Howe_79} from the Rayleigh conductivity of a thin circular aperture.   We have thus shown that in the study of a bias flow through a tilted perforation carried out by \citet{2007_AIAA_Eldredge}, a better match between the    theoretical model and the numerical simulations could have been obtained if the exact expression of the Rayleigh conductivity had been used. This kind of methodology could also possibly be applied to the issue of surface waves' attenuation in ice-covered ocean, since since it seems clear that the usual method of determining it (linear scattering theory) underestimates the attenuation of longer waves \citep{2012_Bennetts}, requiring the addition of a somewhat artificial viscosity to the model \citep{2012_Bennetts,2006_Wang}. It might be interesting to try a simple model of a floating elastic plate with periodic holes for
brine drainage (perhaps with some internal friction) to see how much attenuation this produces.

Then, in section \ref{section:reseau}, we  performed an asymptotic expansion on the scattering problem with respect to the small ratio $\delta$ between the characteristic sizes of the perforation and of the lattice cell. We derived accurate first-order and second-order expansions of the reflection and transmission coefficients, which are expressed in terms of the Rayleigh conductivity of the perforation. These expressions extend those previously given by \citet{1973_JFM_Leppington} and retrieve those obtained by matched asymptotic expansions in \cite{bendali2012}:  they are valid for thick plates with an arbitrary grating of perforations with various geometries. It should be noted that the only constraint required on the size of the lattice cell is that it has to be smaller than half of the acoustic wavelength, contrary to the conventional long wavelength assumption. 

Moreover, we have shown that it is necessary to modify the $\bigcirc(\delta)$ and $\bigcirc(\delta^2)$ expansions of the reflection and transmission coefficients to ensure the conservation of acoustic energy. Thus, the asymptotics provide a first-order and second-order expression of the effective compliance of the perforated plate. The first-order-accurate compliance is the same as the classical expression, \ie the averaging of the Rayleigh conductivity of a single perforation by the area of the lattice unit cell. The accurate second-order compliance involves corrections by a coefficient that depends on the shape of the lattice. Finally, numerical calculations on a typical plate have shown that the second-order corrections, and therefore the effects of the lattice shape, can be neglected without any significant loss of accuracy.

\section*{Acknowledgements}
 Part of this work was funded by the French National Research Agency (ANR) under grant no. ANR$-08-$SYSC$-001$.


\appendix

\section{Details for computation of the functions involved in the Dirichlet and Kelvin principles when the openong is elliptical in shape}\label{sec:calculs}

This appendix discusses the calculations described in section \ref{sec:tiltedelliptical}, and especially the solutions of problem (\ref{eq:pbpotential}) with 
either a Dirichlet or a Neumann boundary condition on the elliptical hole. The elliptical opening in completely outside the scope of 
the results reported in \cite{laurens2012}. Actually, the method of \citet{copson1947} used in  \cite{laurens2012} can no longer be applied. 
To overcome this difficulty, we use an approach similar to the one introduced in \citet{1973_JFM_Leppington} and reused  in \cite{laurensAML}. 
This requires the use of a spheroidal coordinate system and the decomposition of functions into even Legendre functions. In  \cite{laurensAML}, 
the energy is computed for a specific polynomial of degree $1$ source density, while in this work, we determine the source density from the potential it
 generates.  In this section, we thus present new results concerning the determination of the source density generating a linear potential on a thin disk.

\subsection{Solution of problem (\ref{eq:pbpotential}) with a Dirichlet boundary condition on the elliptical hole}

We are solving the problem (\ref{eq:pbpotential}) applied to the potentials $w_{\epss}$ and $t_{\epss}$, which verify the boundary conditions (\ref{eq:30}). The integrals
\begin{equation}\label{eq:21bis}
I_{f_{\epss}} =  \int_{\pm x_3 > h/2} \big|\bnabla f_{\epss}(\bx) \big| ^2 \dif \bx = \int_\sA \rho_{f_{\epss}}(\bx) \cdot f_{\epss}(\bx) \; \dif \bx',
\end{equation}
 with either $f_{\epss}=w_{\epss}$ or  $f_{\epss}=t_{\epss}$, solution to  (\ref{eq:pbpotential}) and boundary condition (\ref{eq:30}), give  the upper bound of the Rayleigh conductivity. This problem can be solved using the method of \cite{copson1947} for a circular opening. But if the opening is elliptical, this method can no longer be applied.  To determine  $\rho_{w_{\epss}}$ and $\rho_{t_{\epss}}$, we use the same mathematical approach as in \cite{laurensAML}. We consider the spheroidal coordinate system:
\begin{equation}\label{eq:3}
  x_1=a\, \sin\theta' \, \cos\varphi'\: \mbox{ and } \: x_2=b\, \sin\theta'\,  \sin\varphi',  \mbox{ with  }\:  \theta' \in [0,\pi/2], \: \varphi' \in [0,2\pi],
\end{equation}
to obtain a parameterization of $\sA$ in terms of the unit
half-sphere. In that case, there is a block diagonalisation of $1/ |\bx - \by|$ on a well-chosen spectral basis for the half-sphere involving the even 
Legendre functions $Q_n^m$, meaning that any function $\tau$ can be decomposed under the following form:
\begin{equation}
  \label{eq:sphericalharmonics}
  \tau(\theta', \varphi') = \sum_{n=0}^{\infty} \:\: \sum^{|m| \leq n}_{ n-m \: { \rm{even} } }  \tau_n^m Q_n^m(\cos \theta') e^{i m \varphi'}.
\end{equation}

We denote by $g_{\epss}(\theta', \varphi') = \cos \theta' \rho_{f_{\epss}}(\theta', \varphi')$ the modified source density function.  Formula $(9)$ of \cite{laurensAML} states that the coefficients $f_n^m$ and  $g_n^m$
 of the spheroidal expansions (\ref{eq:sphericalharmonics}) of $f_{\epss}$ and $g_{\epss}$ are related through 
\begin{equation}
  \label{eq:fnm}
  f_n^m = \frac{\sqrt{ab}}{2} \:\: \sum^{|m'| \leq n}_{ n-m' \: { \rm{even} } }  d_{mm'}^n g_n^m, \footnote{\scriptsize Formula  $(9)$ of \cite{laurensAML}  is actually $ f_n^m = \frac{\sqrt{ab}}{4} \sum^{|m'| \leq n}_{n-m' \rm{even}}   d_{mm'}^n g_n^m$ because the potential formula $(1)$ of \cite{laurensAML}  is  $f_{\epss}(\bx',x_3) = \frac{1}{4 \pi} \int_\sA \rho_{f_{\epss}(\by)} /|\bx' - \by| \dif \by$. In this paper we have used the same convention as in \cite[formula(2.22)]{1973_JFM_Leppington}, see (\ref{eq:electrostaticpotential}).}
\end{equation}
with the coefficients $d_{mm'}^n$ given by Formula $(6)$ of  \cite{laurensAML}, or equivalently by \begin{equation}
  \label{eq:dmmn}
d_{mm'}^n  =  4 \:  \frac{Q_n^m(0)Q_n^{m'}(0)}{2n+1} \: \sqrt{\frac{b}{a}} \int_0^{\pi/2} \frac{ e^{i (m-m') \varphi} }{\sqrt{ 1 - \epss^2 sin^2 \varphi }}\dif \varphi.
\end{equation}

But unlike in \cite{laurensAML}, we need to compute the source density functions $\rho_{w_{\epss}}$ and $\rho_{t_{\epss}}$. To do so, first,
we fully expand the potential $w_{\epss}$ on $\sA$. Since, $w_{\epss}$ is constant on $\sA$,  its spheroidal expansion (\ref{eq:sphericalharmonics}) is:
\begin{equation}
  \label{eq:f00}
   w_{\epss}(\theta', \varphi') =  f_0^0 \: Q_0^0(\cos \theta').
\end{equation}

Secondly, we develop (\ref{eq:fnm}) by considering only the $0$-order Legendre functions
\begin{equation}
  f_0^0 =  (\sqrt{ab} / 2) \: d_{00}^0 \:  g_0^0.
\end{equation}

Then we solve
\begin{equation}
  w_{\epss}(\theta', \varphi') =   (\sqrt{ab} / 2) \: d_{00}^0 \:  g_0^0 \: Q_0^0(\cos \theta').
\end{equation}

Since $Q_0^0 (\cos \theta') = 1/ \sqrt{2}$, $w_{\epss}(\theta', \varphi') = 1$ and $ d_{00}^0 = 2 \sqrt{b/a} \:K(\epss)$, we have
\begin{equation}
  g_0^0 =  \sqrt{2} / (b \: K(\epss) ).
\end{equation}

Given that $ \cos \theta' = \sqrt{ 1 -  \sin^2 \theta'} = \sqrt{ 1 - x_1^{2}/a^{2}-x_2^{2}/b^{2} }$, it follows that the source density function is 
\begin{equation}
  \label{eq:wr}
  \rho_{w_{\epss}}(x_1, x_2)  = \frac{g_{\epss}(\theta', \varphi')}{ \cos \theta' } =  \frac{  g_0^0 \: Q_0^0(\cos \theta')}{ \cos \theta' }  =  \frac{1} {b K(\epss)}  \frac{1}{ \sqrt{ 1 - x_1^{2}/a^{2}-x_2^{2}/b^{2} }}.
\end{equation}

To compute the source density $  \rho_{t_{\epss}}$ that generates a degree $1$ polynomial potential $t_{\epss}$, we expand the potential thanks to (\ref{eq:sphericalharmonics}), meaning 
\begin{equation}
   \label{eq:f00}
   t_{\epss}(\theta', \varphi') =  f_1^1Q_1^1(\cos \theta')e^{i \varphi'} + f_1^{-1}Q_1^{-1}(\cos \theta')e^{ - i \varphi'}.
\end{equation}

Then, keeping all the $0$ et $1$-order Legendre functions in the expansion (\ref{eq:fnm}) yields
\begin{systeme}
  f_1^1 &  =  (\sqrt{ab} / 2) \: \left(  d_{11}^1 \:  g_1^1 +  d_{1,-1}^1 \:  g_1^{1} \right)\\
 f_1^{-1} & =  (\sqrt{ab} / 2) \: \left(  d_{-1,1}^1 \:  g_1^1 +  d_{1,-1}^1 \:  g_1^{-1} \right).
\end{systeme}

 The coefficients $d_{mm'}^n$ are given by (\ref{eq:fnm}):
 \begin{equation}
   d_{11}^1 = 2 \sqrt{ b / a } \: K(\epss) \quad \mbox{ and }  \quad  d_{1,-1}^1 =  d_{-1,1}^1 = 2 \sqrt{ b / a } \:( K(\epss) - 2 D(\epss) ).
 \end{equation}

Since $ t_{\epss}(\theta', \varphi') = x_1 = a\, \sin \varphi' \, \cos \theta'$ and $Q_1^{-1} (\cos \theta') = -  Q_1^{1} (\cos \theta') =  \sqrt{3} \: \sin \theta' / 2$, we find that
\begin{systeme}  \label{eq:f11}  
& f_1^{-1} = -   f_1^1 =  a / \sqrt{3}.\\
&   f_1^1  =  (b/4)\: \left[ \:  K(\epss)  \: g_1^1 +  ( K(\epss) - 2 D(\epss) ) \: g_1^{-1} \right], \\ 
& f_1^{-1} = (b/4) \: \left[ \: ( K(\epss) - 2 D(\epss) )  \: g_1^1 +  K(\epss)  \: g_1^{-1} \right].
\end{systeme}

As  $f_1^{-1}+   f_1^1 =0$, it comes that $g_1^{-1} +  g_1^1 = 0$ and
\begin{equation}
   g_1^{-1}= - g_1^1= - a / (  \sqrt{3} b\: D(\epss)).
\end{equation}

Therefore,  the source density function $\rho_{t_{\epss}}$ that generated a potential which is a multiple of $x_1$ on an elliptic disk is
\begin{equation}
  \label{eq:tr}
   \rho_{t_{\epss}}( x_1, x_2)  =  \frac{1} {b D(\epss)}  \frac{x_1}{ \sqrt{1 - x_1^{2}/a^{2}-x_2^{2}/b^{2} }}.
\end{equation}

The limiting case  $\epss = 0$, with then $a=b=r$, is in agreement with the results obtained by \cite{copson1947} for a circle.

With the analytical expressions of the source density $\rho_{w_{\epss}}$ and
$\rho_{t_{\epss}}$, we can compute the following integrals  from (\ref{eq:21bis}) that have to be applied for the Kelvin and Dirichlet principles:
\begin{equation}\label{eq:integrals}
 \int_{\pm x_3>h/2} \big|\bnabla w_{\epss}(\bx) \big| ^2\dif \bx= \pi a \frac{K(0)}{K(\epss)} \quad  \mbox{ and }  \: 
 \int_{\pm x_3>h/2} \big| \bnabla t_{\epss}(\bx) \big|^2\;\dif \bx=\frac{8}{3} a^3 \frac{D(0)}{D(\epss)}.
\end{equation}

\subsection{Solution of problem (\ref{eq:pbpotential}) with a Neumann boundary condition on the elliptical hole}

We consider the problem (\ref{eq:pbpotential}) with boundary conditions (\ref{eq:25}). The corresponding potential  $z_{\epss}$ is associated to a constant source density function $\rho_{z_{\epss}} = 1/2$ (see Sec. 3.4.1.3, p.144 of \citet{sauterschwab}).  In this case, it has been shown by \cite{laurensAML} that:
\begin{equation}\label{eq:25bis}
\int_{\pm x_3> h/2} \big| \bnabla z_{\epss}(\bx) \big|^2\;\dif \bx=\frac{8}{3} ab^2  \frac{K(\epss)}{K(0)}.
\end{equation}

This result is obtained by using theorem $1.1$  with $\alpha_0 = 1/2, \alpha_1=\alpha_2 = 0$ and the potential defined by (\ref{eq:electrostaticpotential}). As already enlighted in the footnote (\ref{eq:fnm}), the potential formula differs from the one used in \cite{laurensAML} and thus there is a factor 2 of difference between expression (\ref{eq:25bis}) and the one given in \cite{laurensAML}.


\bibliographystyle{jfm}
\bibliography{plaque-perforee}

\end{document}